\documentclass[12pt, reqno]{amsart}
\usepackage{amssymb, latexsym}

\newtheorem{thm}{Theorem}[section]
\newtheorem{corollary}[thm]{Corollary}
\newtheorem{lemma}[thm]{Lemma}
\newtheorem{theorem}[thm]{Theorem}
\newtheorem{definition}[thm]{Definition}
\newtheorem{prop}[thm]{Proposition}

\numberwithin{equation}{section}

\newcommand{\al}{\alpha}
\renewcommand{\b}{\beta}
\renewcommand{\c}{\gamma}
\newcommand{\de}{\delta}
\newcommand{\e}{\varepsilon}
\newcommand{\la}{\lambda}
\renewcommand{\phi}{\varphi}

\renewcommand{\d}{\partial}
\newcommand{\R}{{\mathbb R}}

\newcommand{\br}[1]{\left\langle #1 \right\rangle}
\newcommand{\Case}[1]{\noindent \underline{Case #1:}}
\newcommand{\Step}[1]{\noindent \underline{Step #1:}}
\newcommand{\Subcase}[1]{\noindent \underline{Subcase #1:}}

\renewcommand{\qed}{\rule{3mm}{3mm}}
\renewenvironment{proof}
    {\vspace{1mm}\noindent{\bf Proof.}}
    {\hspace*{\fill} $\qed$\vspace{1mm}}
\newenvironment{proof_of}[1]
    {\vspace{1mm}\noindent {\bf Proof of #1.}}
    {\hspace*{\fill} $\qed$\vspace{1mm}}

\begin{document}
\title[Weighted Decay Estimates for the Wave Equation]{Weighted Decay Estimates for the Wave Equation with Radially
Symmetric Data}
\author{Paschalis Karageorgis}
\address{Brown University, Box 1917, Providence, RI 02912}
\email{petekara@math.brown.edu}



\begin{abstract}
We study the homogeneous wave equation with radially symmetric data in $n\geq 4$ space dimensions.  Using some new
integral representations for the Riemann operator, we estimate the $L\sp \infty$-norm of the solution.  Our results
refine those of Kubo \cite{Kb1, Kb2} in odd space dimensions as well as those of Kubo and Kubota \cite{KKe} in even space
dimensions. However, our approach does not really depend on the parity of $n$.
\end{abstract}
\maketitle

\section{Introduction}

In this paper, we study the homogeneous wave equation
\begin{equation}\label{we}
\left\{
\begin{array}{rll}
\d_t^2 u_0 - \Delta u_0 &= 0 \quad\quad \text{in\: $\R^n \times (0,\infty)$} \\
u_0(x,0) &= \phi(x) \\
\d_t u_0(x,0) &= \psi(x)
\end{array}\right.
\end{equation}
with radially symmetric data in $n$ space dimensions.  Our goal is to obtain \textit{a priori} decay estimates for the
solution when generic assumptions are imposed on the initial data.  Such estimates have been used by several authors in
the existence theory for the semilinear wave equation $\d_t^2 u - \Delta u = |u|^p$, where $p>1$.  In a later paper, we
are going to utilize our results to address the existence of global solutions to the more general nonlinear wave equation
with potential $\d_t^2 u - \Delta u = |u|^p - V(x)\cdot u$, where $p>1$ and $V(x)$ is radially symmetric.

First, consider the solution $u_0$ of \eqref{we} when the initial data are such that
\begin{equation*}
\sum_{|\al|\leq 3} \:|\d_x^\al \phi(x)| + \sum_{|\al|\leq 2} \:|\d_x^\al \psi(x)| \leq \e (1+|x|)^{-k-1}
\end{equation*}
for some $\e>0$ and $k\geq 0$.  In low space dimensions $n=2,3$, sharp decay estimates for $u_0$ were obtained by Kubota
\cite{Ku} and independently by Tsutaya \cite{Ts1, Ts3}; see also the earlier work of Asakura \cite{As}. In these papers,
the assumption of radial symmetry was not needed to control the $L^\infty$-norm of $u_0$.  To obtain similar decay
estimates in higher dimensions $n\geq 4$, however, one has to consider initial data that are either more regular or else
radially symmetric. In what follows, we focus on the latter case and allow a generic regularity assumption.  Thus, our
goal is to study the homogeneous problem
\begin{equation}\label{he}
\left\{
\begin{array}{rll}
\d_t^2 u_0 - \d_r^2 u_0  - \dfrac{n-1}{r} \cdot \d_r u_0 &= 0 \quad\quad \text{in\: $\R_+^2= (0,\infty)^2$} \\
u_0(r,0) &= \phi(r) \\
\d_t u_0(r,0) &= \psi(r)
\end{array}\right.
\end{equation}
for a fixed integer $n\geq 4$.

Before we state our main result, however, let us first introduce some notation.  Given an integer $n\geq 4$, we define
the parameters $a,m$ according to the formula
\begin{equation}\label{am}
(a,m) = \left\{
\begin{array}{ccl}
\left( 1\:,\: \frac{n-3}{2} \right) &&\text{if\, $n$ is odd} \\
\left( \frac{1}{2} \:,\: \frac{n-2}{2} \right) &&\text{if\, $n$ is even.}
\end{array}\right.
\end{equation}
Note that $m\geq 1$ whenever $n\geq 4$ and that the sum $a+m= (n-1)/2$ is independent of the parity of $n$.  Also, we
shall frequently use the bracket notation $\br{\la} = 1+|\la|$ for each $\la\in \R$.  Now, consider the solution to
\eqref{he} when the initial data are such that
\begin{equation}\label{data}
\sum_{s=0}^{l+1} \la^s \,|\phi^{(s)}(\la)| + \sum_{s=0}^l \la^{s+1} \,|\psi^{(s)}(\la)| \leq \e\la^{l-m} \br{\la}^{m-l-k}
\end{equation}
for some $\e>0$, $k\geq 0$ and some integer $1\leq l\leq m$.  Note that this condition allows the data to be singular at
the origin and automatically holds in the case that
\begin{equation*}
\sum_{s=0}^{l+1} \br{\la}^s \,|\phi^{(s)}(\la)| + \sum_{s=0}^l \br{\la}^{s+1} \,|\psi^{(s)}(\la)| \leq \e \br{\la}^{-k}.
\end{equation*}

Our main result in this paper is the following

\begin{theorem}\label{dec}
Let $n\geq 4$ be an integer and define $a,m$ by \eqref{am}.  Fix an integer $1\leq l\leq m$ and consider functions
$\phi\in \mathcal{C}^{l+1}(\R_+)$ and $\psi\in \mathcal{C}^l(\R_+)$ which are subject to \eqref{data}.  Then the
homogeneous equation \eqref{he} admits a unique solution $u_0\in \mathcal{C}^l(\R_+^2)$ that satisfies the following
estimates when $D= (\d_r, \d_t)$ and $\b$ is any multi-index with $|\b|\leq l$.

\begin{itemize}
\item[(a)]
When $0\leq k < m+a= (n-1)/2$,
\begin{equation*}
|D^\b u_0(r,t)| \leq C_0\e r^{l-|\b|-m} \cdot \br{t-r}^{-|\b|} \br{t+r}^{|\b|-l+m-k}.
\end{equation*}

\item[(b)]
When $k= m+a= (n-1)/2$,
\begin{equation*}
|D^\b u_0(r,t)| \leq C_0\e r^{l-|\b|-m} \cdot \br{t-r}^{-|\b|} \br{t+r}^{|\b|-l-a} \left(1+ \ln \frac{\br{t+r}}{\br{t-r}}
\right).
\end{equation*}

\item[(c)]
When $m+a < k < 2(m+a)= n-1$,
\begin{equation*}
|D^\b u_0(r,t)| \leq C_0\e r^{l-|\b|-m} \cdot \br{t-r}^{m+a-k-|\b|} \br{t+r}^{|\b|-l-a}.
\end{equation*}

\item[(d)]
When $k = 2(m+a)=n-1$,
\begin{equation*}
|D^\b u_0(r,t)| \leq C_0\e r^{l-|\b|-m} \cdot \br{t-r}^{-m-a-|\b|} \br{t+r}^{|\b|-l-a} \cdot (1+ \ln \br{t-r}).
\end{equation*}

\item[(e)]
When $k > 2(m+a)=n-1$,
\begin{equation*}
|D^\b u_0(r,t)| \leq C_0\e r^{l-|\b|-m} \cdot \br{t-r}^{-m-a-|\b|} \br{t+r}^{|\b|-l-a}.
\end{equation*}
\end{itemize}
Besides, the constant $C_0$ that appears above depends solely on $k$ and $n$.
\end{theorem}

When it comes to the decay rates $k< (n-1)/2$ of part (a), our method leads to slightly sharper conclusions which we
establish separately, in Corollary \ref{imp}.  Also, when $n$ is odd, our conclusions may be further improved due to the
strong form of Huygens' principle. Although we shall not bother to prove this explicitly, one has $u_0\in
\mathcal{C}^{l+1}(\R_+^2)$ for odd values of $n$, while the estimate of part (c) holds for any decay rate $k> (n-1)/2$
whatsoever.

Let us briefly describe the role of Theorem \ref{dec} in the existence theory for the nonlinear wave equation $\d_t^2 u -
\Delta u = |u|^p$, where $p>1$.  Since the nonlinearity $|u|^p$ is a $\mathcal{C}^1$ function, we require the same
regularity on the initial data.  Then the special case $l=1$ of our theorem asserts that the solution $u_0$ of the
homogeneous equation is $O(r^{1-m})$ as $r\to 0$.  To control this singularity at the origin, one is forced to impose an
upper bound on $p$, as in \cite{Kb1, KKo, KKe}.  Our theorem actually refines the estimates of \cite{Kb1, KKo, KKe} in
any space dimension $n\geq 4$, however, one may avoid such precise estimates by utilizing the method-driven upper bound
on $p$. Our main contribution here is that, unlike \cite{KKe}, we are placing no restrictions on the decay rate $k$ of
the initial data in even space dimensions.

Similarly, let us consider the nonlinear wave equation $\d_t^2 u - \Delta u= |u|^p$ when $p\geq m$.  If we require the
same regularity on the initial data, the special case $l=m$ of our theorem asserts that the solution $u_0$ of the
homogeneous equation is not singular at the origin.  In this case, the upper bound on $p$ is redundant, but the estimates
of Corollary \ref{imp} are almost necessary for an iteration argument to go through.  Although the slightly weaker
estimates of \cite{Kb2} do suffice, those were only obtained for odd space dimensions.

The interesting feature in our method is that our approach is the same regardless of the parity of $n$.  The plausibility
of such an approach is not evident from previous considerations which depended on various representations of the Riemann
operator for the wave equation.  Here, some new representations will be established to facilitate this kind of an
approach.

\begin{definition}\label{pol}
Given an integer $m\geq 1$, define the $m$th Legendre polynomial by
\begin{equation}\label{Pm}
P_m(x) = \frac{1}{2^m\,m!}\cdot \frac{d^m}{dx^m} \:(x^2-1)^m
\end{equation}
and the $m$th Tchebyshev polynomial by
\begin{equation}\label{Tm}
T_m(x) = \frac{(-1)^m}{(2m-1)!!}\cdot\sqrt{1-x^2} \cdot \frac{d^m}{dx^m}\: (1-x^2)^{m-1/2}.
\end{equation}
\end{definition}

\begin{lemma}[The Riemann operator]\label{hs}
Letting $z(\la,r,t)$ be the rational function
\begin{equation}\label{z}
z(\la,r,t) = \frac{\la^2+ r^2-t^2}{2r\la} \:,
\end{equation}
we define the Riemann operator $L$ as follows.  When $n$ is odd, we set
\begin{equation}\label{Lo}
[Lf](r,t) = \frac{1}{2r^{(n-1)/2}} \int_{|t-r|}^{t+r} \la^{(n-1)/2} f(\la) \cdot P_m(z(\la,r,t))\:d\la
\end{equation}
with $m=(n-3)/2$.  When $n$ is even, on the other hand, we set
\begin{align}\label{Le}
[Lf](r,t) &= \frac{\sqrt2}{2\pi r^{(n-1)/2}} \int_{|t-r|}^{t+r} \int_{z(\la,r,t)}^1 \frac{\la^{(n-1)/2}
f(\la)}{\sqrt{\sigma-
z(\la,r,t)}} \cdot \frac{T_m(\sigma)}{\sqrt{1- \sigma^2}} \:\:d\sigma\,d\la \notag\\
&\quad + \frac{\sqrt2}{2\pi r^{(n-1)/2}} \int_0^{\max (t-r,0)} \int_{-1}^1 \frac{\la^{(n-1)/2} f(\la)}{\sqrt{\sigma-
z(\la,r,t)}} \cdot \frac{T_m(\sigma)}{\sqrt{1- \sigma^2}} \:\:d\sigma\,d\la
\end{align}
with $m= (n-2)/2$.  A solution to the Cauchy problem \eqref{he} is then provided by the formula
\begin{equation}\label{u0}
u_0(r,t) = [L\psi](r,t) + \d_t [L\phi](r,t).
\end{equation}
When $\psi\in \mathcal{C}^l(\R_+)$ and $\phi\in \mathcal{C}^{l+1}(\R_+)$ for some integer $l\geq 1$, this solution
belongs to $\mathcal{C}^l(\R_+^2)$.
\end{lemma}

\begin{proof}
Our assertion that \eqref{u0} defines a $\mathcal{C}^l$ function will be established through the proof of Theorem
\ref{dec}, where the derivatives of $u_0$ will be estimated.  Thus, we need only establish the integral representations
\eqref{Lo} and \eqref{Le} for the Riemann operator.  One may derive those using the equivalent representations of Lamb
\cite{La} or Rammaha \cite{Ra}; see also \cite{Ta2}.  For the sake of completeness, however, we shall include their
derivation.

The odd-dimensional representation \eqref{Lo} is precisely the one that appears in \cite{Ra}.  To prove the
even-dimensional representation \eqref{Le}, we manipulate the formula
\begin{equation*}
[Lf](r,t) = \frac{1}{\pi r^{n/2}} \int_0^t \frac{\rho}{\sqrt{t^2-\rho^2}} \:\int_{|\rho-r|}^{\rho+r} \la^{(n-2)/2} f(\la)
\cdot \frac{T_m(z(\la,r,\rho))}{\sqrt{1-z(\la,r,\rho)^2}} \:\:d\la\,d\rho
\end{equation*}
which appears as equation (6b) in \cite{Ra}.  Here, $m= (n-2)/2$ is an integer and $T_m$ is the $m$th Tchebyshev
polynomial \eqref{Tm}.  Switching the order of integration, one arrives at
\begin{align*}
Lf &= \frac{1}{\pi r^{n/2}} \int_{|t-r|}^{t+r} \int_{|\la-r|}^t \frac{\rho}{\sqrt{t^2-\rho^2}} \cdot \la^{(n-2)/2}
f(\la) \cdot \frac{T_m(z(\la,r,\rho))}{\sqrt{1-z(\la,r,\rho)^2}} \:\:d\rho\,d\la \notag\\
&\quad + \frac{1}{\pi r^{n/2}} \int_0^{\max (t-r,0)} \int_{|\la-r|}^{\la+r} \frac{\rho}{\sqrt{t^2-\rho^2}} \cdot
\la^{(n-2)/2} f(\la) \cdot \frac{T_m(z(\la,r,\rho))}{\sqrt{1-z(\la,r,\rho)^2}} \:\:d\rho\,d\la.
\end{align*}
Note that $z(\la,r,|\la\pm r|)= \mp 1$, $\d_\rho z(\la,r,\rho)= -\rho/(r\la)$ and also
\begin{equation*}
t^2-\rho^2 = 2r\la\cdot \Bigl( z(\la,r,\rho) - z(\la,r,t) \Bigr).
\end{equation*}
Once we now use the substitution $\sigma = z(\la,r,\rho)$ in the integrals above, we obtain \eqref{Le}.
\end{proof}

In the remaining of this paper, we shall proceed as follows. In section \ref{be}, we establish some basic facts about the
various functions that appear in the previous lemma. In section \ref{Rm}, we combine these facts to estimate the Riemann
operator and its derivatives. Finally, section \ref{fr} is devoted to the proof of Theorem \ref{dec} and its improved
version, Corollary \ref{imp}.

\section{Basic Estimates and Facts}\label{be}
Our main goal in this section is to collect a few basic facts about the Riemann operator of Lemma \ref{hs}.  Its integral
representation involves the rational function \eqref{z} regardless of the parity of $n$, so we intend to focus on this
function first. In our next lemma, we estimate its derivatives using a rather painstaking approach.  However, we do need
such an approach in order to gain certain cancellations in our subsequent treatise of the Riemann operator.

\begin{lemma}\label{dz}
With $(r,t)\in \R_+^2$ arbitrary and $z= z(\la,r,t)$ as in \eqref{z}, the sharp estimate
\begin{equation}\label{dz1}
|\d_\la^i z| \leq C(i) \cdot \frac{\la^2 + |t^2-r^2|}{r\la^{i+1}}\:, \quad\quad \la\geq 0
\end{equation}
holds for each integer $i\geq 0$, and the general estimate
\begin{equation}\label{dz2}
|\d_\la^i \d_r^j \d_t^k z| \leq C(i,j,k) \cdot \la^{-1-i} \,r^{-1-j} \cdot (t+r)^{2-k}, \quad\quad 0\leq \la\leq t+r
\end{equation}
holds for all integers $i,j,k\geq 0$.  In the case that $t\geq r$, one also has the auxiliary estimate
\begin{equation}\label{dz3}
|D^\b z| \leq C(\b) \cdot \frac{\max(r^2 + |t^2-\la^2|, rt)}{\la r^{|\b|+1}}\:, \quad\quad 0\leq \la\leq t+r
\end{equation}
with $D= (\d_r, \d_t)$ and $\b$ an arbitrary multi-index.
\end{lemma}

\begin{proof}
Our first assertion \eqref{dz1} follows trivially once we explicitly compute
\begin{equation}\label{dzl}
z = \frac{\la^2 -t^2 +r^2}{2r\la} \:\:;\:\: \d_\la z= \frac{\la^2 +t^2 - r^2}{2r\la^2} \:\:;\:\: \d_\la^i z = C(i) \cdot
\frac{t^2- r^2}{r\la^{i+1}} \quad \text{if\, $i\geq 2$}.
\end{equation}
Let us now turn to our second assertion \eqref{dz2}.  According to Leibniz' rule, we have
\begin{equation*}
|\d_\la^i \d_r^j \d_t^k z| \leq C \sum_{i_1=0}^i \sum_{j_1=0}^j \la^{-1-i_1} \,r^{-1-j_1} \cdot |\d_\la^{i-i_1}
\d_r^{j-j_1} \d_t^k (\la^2-t^2+r^2)|.
\end{equation*}
Since $\la,r,t\leq t+r$ by assumption, we may thus deduce the desired estimate
\begin{align*}
|\d_\la^i \d_r^j \d_t^k z|
&\leq C \sum_{i_1=0}^i \sum_{j_1=0}^j \la^{-1-i_1} \,r^{-1-j_1} \cdot (t+r)^{2-(i-i_1)-(j-j_1)-k} \\
&\leq C\la^{-1-i} \,r^{-1-j} \cdot (t+r)^{2-k}.
\end{align*}
Finally, we prove our last assertion \eqref{dz3}. In the case that $\b= (j,0)$, the inequality
\begin{equation*}
|D^\b z| = |\d_r^j z| \leq C(j) \cdot \frac{r^2 + |t^2-\la^2|}{\la r^{j+1}}\:, \quad\quad \la\geq 0
\end{equation*}
follows by \eqref{dz1} because $z$ remains unchanged when the roles of $\la$ and $r$ are interchanged.  To settle the
remaining case $\b= (j,k)$ with $k\geq 1$, we resort to our general estimate \eqref{dz2}.  Since $t\geq r$ by assumption,
$t+r$ is equivalent to $t$, so we get
\begin{equation*}
|D^\b z| = |\d_r^j \d_t^k z| \leq C\la^{-1} \,r^{-1-j} \cdot t^{2-k} \leq C\la^{-1} \,r^{-j-k} \cdot t
\end{equation*}
because $t\geq r$ and $k\geq 1$.  This also completes the proof of our last assertion \eqref{dz3}.
\end{proof}

\begin{corollary}\label{dzb}
With $D_*= (\d_\la, \d_r, \d_t)$ and $\al$ any multi-index, one has
\begin{equation*}
|D_*^\al z(\la,r,t)| \leq C(\al) \cdot \la^{-|\al|}
\end{equation*}
whenever $0< t \leq 2r$ and $|t-r| \leq \la \leq t+r$.
\end{corollary}

\begin{proof}
In the case that $\al= (i,j,k)$ with $j+k\geq 1$, our general estimate \eqref{dz2} gives
\begin{equation*}
|D_*^\al z(\la,r,t)| \leq C\la^{-1-i} \,r^{-1-j} \cdot (t+r)^{2-k} \leq C\la^{-1-i} \,r^{1-j-k}\leq C\la^{-i-j-k}
\end{equation*}
because $r\leq t+r\leq 3r$ and $\la \leq t+r \leq 3r$ by assumption.

In the case that $\al= (i, 0, 0)$, on the other hand, our sharp estimate \eqref{dz1} gives
\begin{equation*}
|D_*^\al z(\la,r,t)| = |\d_\la^i z(\la,r,t)| \leq C(i) \cdot \frac{\la^2 + |t^2-r^2|}{r\la^{i+1}} \leq C(i) \cdot
\frac{\la + t+r}{r\la^i} \leq C\la^{-i}
\end{equation*}
since $|t-r| \leq \la$ and $\la+t+r \leq 2(t+r) \leq 6r$.  In either case then, the result follows.
\end{proof}

In our next corollary, we concern ourselves with the derivatives of $1/(\d_\la z)$.  This function will arise as soon as
we integrate by parts the integrals of Lemma \ref{hs}, namely the ones that appear in the explicit representation of the
Riemann operator. To estimate its derivatives in a rather precise manner, we first recall \eqref{dz1} and \eqref{dzl},
according to which
\begin{equation*}
(\d_\la z)^{-1} = \frac{2r\la^2}{\la^2 +t^2 - r^2}\:,\quad\quad |\d_\la^i z| \leq C(i)\cdot \frac{\la^2 +
|t^2-r^2|}{r\la^{i+1}}\:.
\end{equation*}
Keeping such expressions intact, one does not have to estimate them separately while dealing with their product. This is
also the main idea in the following proof, where special emphasis is laid on the derivatives with respect to $\la$.

\begin{corollary}\label{dlz}
With $z= z(\la,r,t)$ as in \eqref{z}, one has
\begin{equation}\label{dlz1}
\left| \d_\la^i \d_r^j \d_t^k (\d_\la z)^{-1} \right| \leq C(i,j,k)\cdot \frac{\la^{2-i} \,r^{1-j}
\,t^{j-1}}{(t-r)^{1+j+k}}
\end{equation}
whenever $0< r\leq t$ and $0\leq \la \leq t+r$.
\end{corollary}

\begin{proof}
Let us set $D= (\d_r, \d_t)$ and $\b= (j,k)$ for ease of notation. By repeated applications of the chain rule, one
obtains an identity of the form
\begin{equation}\label{cr}
D^\b (\d_\la z)^{-1} = \sum_{l=1}^{|\b|} \:(\d_\la z)^{-1-l} \sum_{\b_1+ \ldots + \b_l= \b} C(\b_1,\ldots,\b_l, l) \cdot
\prod_{q=1}^l \left( D^{\b_q} \d_\la z \right),
\end{equation}
where each $\b_q$ is a multi-index.  Further differentiating with respect to $\la$, one then obtains
\begin{equation}\label{df1}
\left| \d_\la^i D^\b (\d_\la z)^{-1} \right| \leq C\sum_{l=1}^{|\b|} \sum_{s=0}^i \left| \d_\la^{i-s} (\d_\la z)^{-1-l}
\right| \sum_{\b_1+ \ldots + \b_l= \b} \left| \d_\la^s \prod_{q=1}^l \left( D^{\b_q} \d_\la z \right) \right|
\end{equation}
by Leibniz' rule.  Once we write $\b_q= (j_q, k_q)$ for each $q$, the innermost sum becomes
\begin{equation*}
\sum_{\b_1+ \ldots + \b_l= \b} \left| \d_\la^s \prod_{q=1}^l \left( D^{\b_q} \d_\la z \right) \right| \leq \sum_{j_1+
\ldots + j_l= j} \: \sum_{k_1+ \ldots + k_l= k} \: \sum_{s_1+ \ldots + s_l= s} \:C \prod_{q=1}^l \left| \d_\la^{s_q+1}
\d_r^{j_q} \d_t^{k_q} z \right|.
\end{equation*}
Meanwhile, our general estimate \eqref{dz2} ensures that
\begin{equation*}
\left| \d_\la^{s_q+1} \d_r^{j_q} \d_t^{k_q} z \right| \leq C\la^{-2-s_q} \,r^{-1-j_q} \cdot (t+r)^{2-k_q}, \quad\quad
0\leq \la\leq t+r
\end{equation*}
so we may combine the last two equations to arrive at
\begin{equation*}
\sum_{\b_1+ \ldots + \b_l= \b} \left| \d_\la^s \prod_{q=1}^l \left( D^{\b_q} \d_\la z \right) \right| \leq C\la^{-2l-s}
\,r^{-l-j} \cdot (t+r)^{2l-k}.
\end{equation*}
Next, we insert this fact in \eqref{df1}.  Since $D= (\d_r, \d_t)$ and $\b= (j,k)$ by above, we find
\begin{equation}\label{df2}
\left| \d_\la^i \d_r^j \d_t^k (\d_\la z)^{-1} \right| \leq C\sum_{l=1}^{j+k} \sum_{s=0}^i \left| \d_\la^{i-s} (\d_\la
z)^{-1-l} \right| \cdot \la^{-2l-s} \,r^{-l-j} \cdot (t+r)^{2l-k}.
\end{equation}
To handle the latter derivatives, we will employ our sharp estimate \eqref{dz1} instead of \eqref{dz2}.  First, we use
repeated applications of the chain rule to establish the inequality
\begin{equation*}
\left| \d_\la^{i-s} (\d_\la z)^{-1-l} \right| \leq C\sum_{m=1}^{i-s} |\d_\la z|^{-1-l-m} \sum_{a_1+ \ldots + a_m= i-s}
\:\prod_{q=1}^m |\d_\la^{a_q+1} z|
\end{equation*}
in analogy with \eqref{cr}.  When it comes to the product, \eqref{dz1} and \eqref{dzl} combine to give
\begin{equation*}
\prod_{q=1}^m |\d_\la^{a_q+1} z| \leq C_1\prod_{q=1}^m \left( \frac{\la^2 + t^2-r^2}{r\la^{a_q+2}} \right) = 2^m
C_1|\d_\la z|^m\cdot \prod_{q=1}^m \la^{-a_q},
\end{equation*}
so we get
\begin{equation*}
\left| \d_\la^{i-s} (\d_\la z)^{-1-l} \right| \leq C|\d_\la z|^{-1-l} \cdot \la^{s-i}.
\end{equation*}
Inserting this fact in \eqref{df2}, we now arrive at
\begin{equation*}
\left| \d_\la^i \d_r^j \d_t^k (\d_\la z)^{-1} \right| \leq C\sum_{l=1}^{j+k} |\d_\la z|^{-1-l} \cdot \la^{-2l-i}
\,r^{-l-j} \cdot (t+r)^{2l-k}.
\end{equation*}
According to our computation \eqref{dzl}, we also have
\begin{equation*}
|\d_\la z|^{-1-l} = \left( \frac{2\la^2r}{\la^2+t^2-r^2}\right)^{l+1} \leq C\la^{2l+2} \,r^{l+1} \cdot (t+r)^{-l-1}
(t-r)^{-l-1},
\end{equation*}
whence
\begin{equation*}
\left| \d_\la^i \d_r^j \d_t^k (\d_\la z)^{-1} \right| \leq C\sum_{l=1}^{j+k} \la^{2-i} \,r^{1-j} \cdot (t+r)^{l-j-k+j-1}
(t-r)^{-l-1}.
\end{equation*}
Since $l\leq j+k$ within the last sum, this trivially gives
\begin{equation*}
\left| \d_\la^i \d_r^j \d_t^k (\d_\la z)^{-1} \right| \leq C\la^{2-i} \,r^{1-j} \cdot (t+r)^{j-1} (t-r)^{-1-j-k}.
\end{equation*}
Besides, $t+r$ is equivalent to $t$ whenever $t\geq r$, so the desired estimate \eqref{dlz1} follows.
\end{proof}

The last fact we need to treat the Riemann operator in odd dimensions is also the most crucial one and appears in our
next lemma.  Here, we use the integral representation \eqref{Lo} to derive some new representations \eqref{Lo2}.  The
latter are only valid in the interior of the light cone, however, they are less singular at $r=0$.

\begin{lemma}\label{Lon}
Let $n\geq 5$ be an odd integer and $a,m$ be as in \eqref{am}. Suppose that $f_0\in \mathcal{C}^l(\R_+)$ and $f_1\in
\mathcal{C}^{l+1}(\R_+)$ for some integer $1\leq l\leq m$.  When $t\geq r$, the Riemann operator \eqref{Lo} is then
subject to an identity of the form
\begin{equation}\label{Lo2}
\d_t^i [Lf_i](r,t) = \frac{(-1)^j}{2} \int_{t-r}^{t+r} [H_{ij} f_i](\la,r,t)\cdot r^{-m-a} P_{jm}(z(\la,r,t)) \:d\la
\end{equation}
for $i=0,1$ and all integers $i\leq j\leq l$.  Here, $P_{jm}$ denotes a polynomial \eqref{Pjm} which vanishes with order
$j$ at each of $\pm 1$, while $H_{ij}$ denotes a linear operator which acts on functions of $\la$ and is defined by
either \eqref{H0j} or \eqref{H1j}. If we let $D= (\d_r, \d_t)$ and $\b$ be any multi-index, then we also have the
estimate
\begin{equation}\label{He}
\left| D^\b [H_{ij}f] \right| \leq C(\b,j) \cdot \frac{\la^j \,r^{j-|\b|} \,t^{|\b|-j}}{(t-r)^{i+j+|\b|}} \cdot
\sum_{s=0}^{i+j} \la^{m+a+s} \,|f^{(s)}(\la)|, \quad\quad 0\leq \la\leq t+r
\end{equation}
for $i=0,1$ and any integer $j\geq 0$, provided that $t\geq r$ as above.
\end{lemma}

\begin{proof}
First, we shall use induction on $j$ to prove \eqref{Lo2} for the case $i=0$ and $0\leq j\leq l$. Explicitly, the
polynomials $P_{jm}$ of interest are given by the formula
\begin{equation}\label{Pjm}
P_{jm}(x) \equiv \frac{1}{2^m\,m!}\cdot \frac{d^{m-j}}{dx^{m-j}} \:(x^2-1)^m, \quad\quad 0\leq j\leq m.
\end{equation}
Note that $P_{0m}$ is merely the Legendre polynomial \eqref{Pm}. Besides, $P_{jm}$ vanishes with order $j$ at each of
$\pm 1$, and we also have $P_{jm} = P_{j+1,m}'$ for each $j$.  As for the linear operators $H_{0j}$ we are going to need,
those are given by the formula
\begin{equation}\label{H0j}
[H_{0j} f](\la,r,t) = \left[ \frac{\d}{\d\la} \:\frac{1}{\d_\la z(\la,r,t)} \right]^j \Bigl( \la^{m+a} f(\la) \Bigr),
\end{equation}
where $z(\la,r,t)$ is the rational function \eqref{z}.  We may introduce them for any integer $j\leq l$ and any function
$f\in \mathcal{C}^l(\R_+)$.  Let us now proceed to the induction argument.  Since $t\geq r$ by assumption, the Riemann
operator \eqref{Lo} takes the form
\begin{equation*}
[Lf_0](r,t) = \frac{1}{2} \int_{t-r}^{t+r} \la^{(n-1)/2} f_0(\la) \cdot r^{-(n-1)/2} \,P_{0m}(z(\la,r,t))\:d\la.
\end{equation*}
Moreover, $(n-1)/2= m+a$ by our definition \eqref{am}, so the desired identity
\begin{equation}\label{ind}
[Lf_0](r,t) = \frac{(-1)^j}{2} \int_{t-r}^{t+r} [H_{0j} f_0](\la,r,t) \cdot r^{-m-a} \,P_{jm}(z(\la,r,t))\:d\la
\end{equation}
does hold when $j=0$. Suppose it holds for some $j\geq 0$.  As we have already remarked, the polynomials \eqref{Pjm} are
such that
\begin{equation*}
P_{jm}(z(\la,r,t)) = P_{j+1,m}'(z(\la,r,t)) = \frac{1}{\d_\la z(\la,r,t)} \cdot \d_\la P_{j+1,m}(z(\la,r,t)).
\end{equation*}
Integrating \eqref{ind} by parts, we thus get to replace the leftmost factor in the integrand by
\begin{equation*}
- \frac{\d}{\d\la} \left( \frac{1}{\d_\la z(\la,r,t)} \cdot [H_{0j} f_0](\la,r,t) \right) = -[H_{0,j+1} f_0](\la,r,t)
\end{equation*}
and the polynomial $P_{jm}$ by $P_{j+1,m}$.  This allows us to finish the inductive proof of \eqref{ind} as long as no
boundary terms arise in the process.  In fact, $P_{j+1,m}$ vanishes at each of $\pm 1$ and our definition \eqref{z} gives
$z(t\pm r,r,t)= \pm 1$, hence no boundary terms arise, indeed.

Next, we establish \eqref{Lo2} for the case $i=1$ and $1\leq j\leq l$.  The linear operators $H_{1j}$ we are going to
need are given by the formula
\begin{equation}\label{H1j}
[H_{1j}f](\la,r,t) = \d_t [H_{0j}f](\la,r,t) - \d_\la \Bigl( (\d_\la z)^{-1} \cdot \d_t z\cdot [H_{0j}f](\la,r,t) \Bigr).
\end{equation}
We may introduce them for any integer $j\leq l$ and any function $f\in \mathcal{C}^{l+1}(\R_+)$.  Let us now employ the
identity \eqref{ind} we just proved to write
\begin{equation*}
[Lf_1](r,t) = \frac{(-1)^j}{2} \int_{t-r}^{t+r} [H_{0j} f_1](\la,r,t) \cdot r^{-m-a} \,P_{jm}(z(\la,r,t))\:d\la.
\end{equation*}
Since $j\geq 1$, we have $P_{jm}(z(t\pm r,r,t))= P_{jm}(\pm 1)= 0$ by above, so we find that
\begin{align*}
\d_t[Lf_1](r,t) &= \frac{(-1)^j}{2} \int_{t-r}^{t+r} \d_t [H_{0j}f_1](\la,r,t) \cdot r^{-m-a} P_{jm}(z(\la,r,t))\:d\la \\
&\quad + \frac{(-1)^j}{2} \int_{t-r}^{t+r} \d_t z\cdot [H_{0j}f_1](\la,r,t) \cdot r^{-m-a} P'_{jm}(z(\la,r,t))\:d\la.
\end{align*}
Recalling our definition \eqref{H1j}, we may then integrate the latter integral by parts to get
\begin{equation*}
\d_t[Lf_1](r,t) = \frac{(-1)^j}{2} \int_{t-r}^{t+r} [H_{1j}f_1](\la,r,t) \cdot r^{-m-a} P_{jm}(z(\la,r,t))\:d\la.
\end{equation*}
This is precisely the desired identity \eqref{Lo2} for the case $i=1$ and $1\leq j\leq l$.

Finally, we turn our attention to the estimate \eqref{He}, which asserts that
\begin{equation}\label{di}
\left|\d_r^{j_0} \d_t^{k_0} [H_{ij}f] \right| \leq C(j_0,k_0,j) \cdot \frac{\la^j \,r^{j-j_0-k_0} \,t^{j_0+k_0
-j}}{(t-r)^{i +j+j_0+k_0}} \cdot \sum_{s=0}^{i+j} \la^{m+a+s} \,|f^{(s)}(\la)|
\end{equation}
for $i=0,1$ and each $0\leq \la\leq t+r$.  Allowing the remaining indices to be all arbitrary, our plan is to establish
the more general estimate
\begin{equation}\label{mg}
\left| \d_\la^{i_0} \d_r^{j_0} \d_t^{k_0} [H_{ij}f] \right| \leq C(i_0,j_0,k_0,j) \cdot \frac{\la^{j-i_0} \,r^{j-j_0}
\,t^{j_0-j}}{(t-r)^{i+j+j_0+k_0}} \cdot \sum_{s=0}^{i+i_0+j} \la^{m+a+s} \,|f^{(s)}(\la)|
\end{equation}
for $i=0,1$ and each $0\leq \la\leq t+r$.  When $i_0=0$, this actually improves \eqref{di} by an extra factor $r^{k_0}
t^{-k_0}$, which is less than $1$ since $r\leq t$ by assumption. In what follows, we may thus focus on the derivation of
\eqref{mg}, instead.

To prove \eqref{mg} for the case $i=0$, we need to check that
\begin{equation}\label{ih}
\left| \d_\la^{i_0} \d_r^{j_0} \d_t^{k_0} [H_{0j}f] \right| \leq C(i_0,j_0,k_0,j) \cdot \frac{\la^{j-i_0} \,r^{j-j_0}
\,t^{j_0-j}}{(t-r)^{j+j_0+k_0}} \cdot \sum_{s=0}^{i_0+j} \la^{m+a+s} \,|f^{(s)}(\la)|
\end{equation}
for each $0\leq \la\leq t+r$.  Let us proceed using induction on $j$.  When $j=0$ and $j_0+k_0\geq 1$, our task is
trivial since $H_{00}f = \la^{m+a} f(\la)$ depends only on $\la$.  When $j=0$ and $j_0=k_0=0$, on the other hand, an
application of Leibniz' rule gives
\begin{equation*}
|\d_\la^{i_0} [H_{00}f]| \leq C(i_0) \cdot \sum_{s=0}^{i_0} \la^{m+a+s-i_0} \,|f^{(s)}(\la)|.
\end{equation*}
This proves \eqref{ih} when $j= 0$, so suppose the same estimate holds for some $j\geq 0$. In view of our definition
\eqref{H0j}, we may then write
\begin{equation*}
[H_{0,j+1}f](\la,r,t) = \d_\la \Bigl( (\d_\la z)^{-1} \cdot [H_{0j}f](\la,r,t) \Bigr).
\end{equation*}
The derivatives of $(\d_\la z)^{-1}$ were treated in Corollary \ref{dlz}, while those of $H_{0j}f$ are subject to our
induction hypothesis \eqref{ih}.  Resorting to Leibniz' rule, one may easily use these facts to finish the inductive
proof of \eqref{ih}, so we shall omit the details.

To establish \eqref{mg} for the case $i=1$, we need to check that
\begin{equation}\label{ih2}
\left| \d_\la^{i_0} \d_r^{j_0} \d_t^{k_0} [H_{1j}f] \right| \leq C(i_0,j_0,k_0,j) \cdot \frac{\la^{j-i_0} \,r^{j-j_0}
\,t^{j_0-j}}{(t-r)^{1+j+j_0+k_0}} \cdot \sum_{s=0}^{i_0+1+j} \la^{m+a+s} \,|f^{(s)}(\la)|
\end{equation}
for each $0\leq \la\leq t+r$.  Here, we may proceed directly starting with the definition
\begin{equation*}
[H_{1j}f](\la,r,t) = \d_t [H_{0j}f](\la,r,t) - \d_\la \Bigl( (\d_\la z)^{-1} \cdot \d_t z\cdot [H_{0j}f](\la,r,t) \Bigr)
\end{equation*}
we introduced in \eqref{H1j}.  The derivatives of $H_{0j}f$ are subject to \eqref{ih}, those of $(\d_\la z)^{-1}$ were
treated in Corollary \ref{dlz}, while our general estimate \eqref{dz2} applies for the derivatives of $z$.  In view of
these facts, our last assertion \eqref{ih2} is now easy to deduce using Leibniz' rule.
\end{proof}

Our next step is to establish an even-dimensional analogue of the previous lemma.  Here, our task is similar but harder,
as the representation \eqref{Le} for the Riemann operator is more subtle. Consider, for instance, the function
\begin{equation*}
U_{0m}(\la,r,t) \equiv \int_{z(\la,r,t)}^1 [\sigma-z(\la,r,t)]^{-1/2} \cdot \frac{T_m(\sigma)}{\sqrt{1-\sigma^2}}
\:\:d\sigma, \quad\quad |t-r|\leq \la\leq t+r
\end{equation*}
that appears in \eqref{Le} as an inner integral.  At the endpoints $\la= r\pm t$, this function behaves rather nicely and
actually attains the value $\pi/\sqrt 2$.  At the endpoint $\la= t-r$, however, the integral above happens to diverge.
Namely, our definition \eqref{z} gives $z(t-r,r,t)= -1$ and this creates a singularity that is not integrable near
$\sigma= -1$. It is worth noting that we only have to deal with the troublesome endpoint when $t\geq r$.  To handle this
case, we shall need to introduce and study the functions
\begin{equation*}
U_{jm}(\la,r,t) \equiv \int_{z(\la,r,t)}^1 [\sigma-z(\la,r,t)]^{j-1/2} \cdot \frac{T_m(\sigma)}{\sqrt{1-\sigma^2}}
\:\:d\sigma, \quad\quad |t-r|\leq \la\leq t+r
\end{equation*}
for all integers $0\leq j\leq m$.  These will play a role similar to that of the polynomials $P_{jm}$ we introduced in
odd dimensions \eqref{Pjm}.

\begin{lemma}\label{Ue}
Fix an integer $m\geq 1$ and let $T_m$ be the polynomial \eqref{Tm}.  With $z= z(\la,r,t)$ given by \eqref{z} whenever
$(r,t)\in \R_+^2$, set
\begin{equation}\label{U}
U_{jm}(\la, r, t) \equiv \int_z^1 (\sigma - z)^{j-1/2} \cdot \frac{T_m(\sigma)}{\sqrt{1- \sigma^2}} \:\:d\sigma
\end{equation}
for each $|t-r|\leq \la\leq t+r$ and each integer $0\leq j\leq m$.

\begin{itemize}
\item[(a)]
At the endpoint $\la= t+r$, the function $U_{jm}$ has a zero of order $j$.

\item[(b)]
Letting $D_*= (\d_\la, \d_r, \d_t)$ and $\al$ be any multi-index, one has
\begin{equation}\label{Ue1}
|D_*^\al U_{0m}(\la,r,t)| \leq C(\al,m) \cdot \left( \frac{\la}{r-t+\la} \right)^{1/2+|\al|} \cdot \la^{-|\al|}
\end{equation}
whenever $0< t \leq 2r$ and $|t-r| \leq \la \leq t+r$.

\item[(c)]
Assume that $t\geq r$.  Letting $D= (\d_r, \d_t)$, one then has
\begin{equation}\label{Ue2}
|D^\b U_{jm}(\la,r,t)| \leq C(m) \cdot \left( \frac{\la r}{t+r} \right)^{1/2-|\b|} \cdot \frac{1}{\sqrt{r-t+\la}}
\end{equation}
for each $|t-r|\leq \la \leq t+r$ and each multi-index $\b$ with $|\b|\leq j\leq m$.
\end{itemize}
\end{lemma}

\begin{proof}
For the first two parts, we use the substitution $\nu = (\sigma -z)/(1-z)$ to write
\begin{equation*}
U_{jm}(\la,r,t) = (1-z)^j \int_0^1 \frac{\nu^{j-1/2}}{\sqrt{1-\nu}} \cdot \frac{T_m(\sigma)}{\sqrt{1+\sigma}} \:\:d\nu,
\quad\quad \sigma= \nu + (1-\nu)z.
\end{equation*}
In view of \eqref{z}, we have $z=1$ when $\la=t+r$, so part (a) is clear.  To prove the estimate of part (b), we take
$j=0$ in the last equation and differentiate to find that
\begin{equation}\label{Ud}
D_*^\al U_{0m}(\la,r,t) = \int_0^1 \frac{1}{\sqrt{\nu(1-\nu)}} \cdot D_*^\al \left( \frac{T_m(\sigma)}{\sqrt{1+\sigma}}
\right) \:d\nu, \quad\quad \sigma= \nu + (1-\nu)z.
\end{equation}
When it comes to the rightmost factor in the integrand, we have
\begin{equation*}
D_*^\al \left( \frac{T_m(\sigma)}{\sqrt{1+\sigma}} \right) = \sum_{s=1}^{|\al|} \:\: \sum_{\al_1+ \ldots + \al_s= \al}
C(\al_1, \ldots, \al_s, s) \cdot \frac{d^s}{d\sigma^s} \left( \frac{T_m(\sigma)}{\sqrt{1+\sigma}} \right) \cdot
\prod_{q=1}^s \left( D_*^{\al_q} \sigma \right).
\end{equation*}
Since $\la$, $r$ and $t$ are as in Corollary \ref{dzb}, an estimate of the form $|D_*^{\al_q} z|\leq C(\al_q)\cdot
\la^{-|\al_q|}$ holds for $z$, hence also for $\sigma = \nu + (1-\nu)z$ because $0\leq \nu \leq 1$.  In particular,
$\sigma$ is bounded and
\begin{equation*}
\left| D_*^\al \left( \frac{T_m(\sigma)}{\sqrt{1+\sigma}} \right) \right| \leq C(\al,m) \cdot (1+\sigma)^{-1/2-|\al|}
\cdot \la^{-|\al|}
\end{equation*}
by above.  Inserting this fact in \eqref{Ud}, we then arrive at
\begin{equation*}
\left| D_*^\al U_{0m}(\la,r,t) \right| \leq C\la^{-|\al|} \cdot \int_0^1 \frac{1}{\sqrt{\nu(1-\nu)}} \cdot
(1+\sigma)^{-1/2 -|\al|} \:\:d\nu, \quad\quad \sigma= \nu + (1-\nu)z.
\end{equation*}
It is easy to check that $|z|\leq 1$ for each $|t-r|\leq \la \leq t+r$, whence
\begin{equation*}
1 + \sigma = 1+ \nu + (1-\nu) z = 1+z + (1-z)\nu \geq 1+z
\end{equation*}
within the region of integration.  Meanwhile, the definition \eqref{z} of $z$ ensures that
\begin{equation*}
1+z = \frac{(r+t+\la)\cdot (r-t+\la)}{2r\la} \geq \frac{r-t+\la}{\la}
\end{equation*}
for each $\la \geq |r-t|$.  Once we now combine the last three equations, we obtain \eqref{Ue1}.

Next, we turn to part (c).  Here, we wish to differentiate $U_{jm}$ a total of $|\b|\leq j$ times.  As the integrand in
\eqref{U} vanishes with order $j-1/2$ at the lower limit of integration, we get
\begin{equation}\label{Db1}
D^\b U_{jm}(\la,r,t) = \int_z^1 D^\b (\sigma - z)^{j-1/2} \cdot \frac{T_m(\sigma)}{\sqrt{1- \sigma^2}} \:\:d\sigma.
\end{equation}
For the leftmost factor in the integrand, repeated applications of the chain rule give
\begin{equation*}
|D^\b (\sigma-z)^{j-1/2}| \leq C(\b) \cdot \sum_{s=1}^{|\b|} \:\: \sum_{\b_1+ \ldots + \b_s= \b} (\sigma -z)^{j-s-1/2}
\cdot \prod_{q=1}^s |D^{\b_q} z|.
\end{equation*}
Moreover, one has $|z|\leq 1$ whenever $|t-r|\leq \la\leq t+r$, so one easily finds
\begin{equation*}
(\sigma - z)^{j-s-1/2} \leq (1-z)^{j-s} \cdot (\sigma-z)^{-1/2} \leq 2^{j-s} \cdot (\sigma-z)^{-1/2}
\end{equation*}
whenever $z\leq \sigma\leq 1$ and $s\leq |\b|\leq j$.  Combining the last two equations, we then obtain
\begin{equation}\label{Db2}
|D^\b (\sigma-z)^{j-1/2}| \leq C\sum_{s=1}^{|\b|} \:\: \sum_{\b_1+ \ldots + \b_s= \b} (\sigma -z)^{-1/2} \cdot
\prod_{q=1}^s |D^{\b_q} z|.
\end{equation}
Since $t\geq r$ by assumption, our auxiliary estimate \eqref{dz3} is applicable here.  In fact, we have
\begin{equation*}
r^2 + |t^2-\la^2| \leq r(r+ t+\la) \leq 2r(r+t), \quad\quad |t-r|\leq \la\leq t+r
\end{equation*}
because $|t-\la|\leq r$ for such $\la$, hence \eqref{dz3} ensures that
\begin{equation}\label{dz5}
|D^{\b_q} z| \leq C(\b_q)\cdot \frac{r(t+r)}{\la r^{|\b_q|+1}} = C(\b_q)\cdot \frac{t+r}{\la r^{|\b_q|}} \:,\quad\quad
|t-r|\leq \la\leq t+r.
\end{equation}
Once we now combine the last equation with \eqref{Db2}, we get
\begin{equation*}
|D^\b (\sigma-z)^{j-1/2}| \leq Cr^{-|\b|} \: \sum_{s=1}^{|\b|} (\sigma -z)^{-1/2} \left( \frac{t+r}{\la} \right)^s \leq
C(\sigma -z)^{-1/2} \left( \frac{t+r}{\la r} \right)^{|\b|}
\end{equation*}
because $\la\leq t+r$.  Inserting this fact in \eqref{Db1}, we thus arrive at
\begin{equation}\label{Ud2}
|D^\b U_{jm}(\la,r,t)| \leq C\left( \frac{t+r}{\la r} \right)^{|\b|} \cdot \int_z^1 (\sigma - z)^{-1/2} \cdot
(1-\sigma^2)^{-1/2} \cdot |T_m(\sigma)| \:d\sigma.
\end{equation}
Note that the polynomial factor is bounded, as $|z|\leq 1$ by above.  If we also neglect a factor of $(1+\sigma)^{-1/2}$
from the integrand, the resulting integral becomes
\begin{equation*}
\int_z^1 (\sigma-z)^{-1/2} \cdot (1-\sigma)^{-1/2} \:d\sigma = \int_0^1 \nu^{-1/2} (1-\nu)^{-1/2} \:d\nu = \pi
\end{equation*}
by means of the substitution $\nu = (\sigma -z)/(1-z)$.  In particular, \eqref{Ud2} leads to
\begin{equation*}
|D^\b U_{jm}(\la,r,t)| \leq C\left( \frac{t+r}{\la r} \right)^{|\b|} \cdot (1+z)^{-1/2}.
\end{equation*}
Besides, the definition \eqref{z} of $z$ is such that
\begin{equation*}
1 + z = \frac{(r+t+\la)\cdot (r-t+\la)}{2\la r} \geq \left( \frac{t+r}{2\la r} \right) \cdot (r-t+\la)
\end{equation*}
for each $\la\geq |t-r|$, so we may combine the last two equations to finally deduce \eqref{Ue2}.
\end{proof}

The functions of the previous lemma are related to the first component of the Riemann operator in even dimensions
\eqref{Le}. To treat the second component, we shall need to study some similar functions that we introduce below
\eqref{W}.  Although their definition resembles the one we had before \eqref{U}, we are now interested in the values
$0\leq \la\leq t-r$.  For such values, the rational function \eqref{z} is no longer bounded and a much more delicate
approach is needed.  In fact, it is only here that the Tchebyshev polynomial \eqref{Tm} becomes important.

\begin{lemma}\label{Wet}
Fix an integer $m\geq 1$ and let $T_m$ be the polynomial \eqref{Tm}.  With $z= z(\la,r,t)$ given by \eqref{z} whenever
$(r,t)\in \R_+^2$, set
\begin{equation}\label{W}
W_{im}(\la, r, t) \equiv \int_{-1}^1 (\sigma - z)^{i-1/2} \cdot \frac{T_m(\sigma)}{\sqrt{1- \sigma^2}} \:\:d\sigma
\end{equation}
for each $0\leq \la\leq t-r$ and each integer $0\leq i\leq m$.

\begin{itemize}
\item[(a)]
At the endpoint $\la= t-r$, the function $W_{im}$ agrees with the function $U_{im}$ of \eqref{U} up to any derivative of
order $i$.

\item[(b)]
Assuming that $0\leq \la\leq t-r$, the estimate
\begin{equation}\label{We1}
|W_{im}(\la,r,t)| \leq C(m) \cdot \left( \frac{r}{t+r} \cdot \frac{\la}{t-r-\la} \right)^{i_0-i+1/2}
\end{equation}
holds for all integers $0\leq i\leq i_0\leq m$.
\end{itemize}
\end{lemma}

\begin{proof}
The only difference between the definitions of $U_{im}$ and $W_{im}$ is that the lower limit of integration is
$z(\la,r,t)$ for the former and $-1$ for the latter.  By \eqref{z}, these two quantities agree when $\la= t-r$, so part
(a) follows easily.  To prove the estimate of part (b), we use the definition \eqref{Tm} of the Tchebyshev polynomial
$T_m$ to first write
\begin{equation*}
W_{im}(\la,r,t) =  C(m) \int_{-1}^1 (\sigma - z)^{i-1/2} \cdot \frac{d^{i_0}}{d\sigma^{i_0}} \left[ \frac{d^{m-i_0}}{d
\sigma^{m -i_0}}\: (1-\sigma^2)^{m-1/2} \right] \:d\sigma.
\end{equation*}
In the case that $i_0\geq 1$, the expression in square brackets vanishes at $\pm 1$.  This means that we may integrate by
parts $i_0$ times to obtain
\begin{equation*}
W_{im}(\la,r,t) =  C(i_0,i,m) \int_{-1}^1 (\sigma - z)^{i-i_0-1/2} \cdot \left[ \frac{d^{m-i_0}}{d \sigma^{m -i_0}}\:
(1-\sigma^2)^{m-1/2} \right] \:d\sigma.
\end{equation*}
Here, the expression in square brackets is integrable on $[-1,1]$ whenever $i_0\geq 0$, so we get
\begin{equation*}
|W_{im}(\la,r,t)| \leq  C(-1 - z)^{i-i_0-1/2}
\end{equation*}
for all integers $0\leq i\leq i_0$.  Since the definition \eqref{z} of $z$ provides the inequality
\begin{equation*}
-1 - z = \frac{r+t+\la}{2r} \cdot \frac{t-r-\la}{\la} \geq \frac{t+r}{2r} \cdot \frac{t-r-\la}{\la}
\end{equation*}
for each $0\leq \la\leq t-r$, we may then combine the last two equations to deduce \eqref{We1}.
\end{proof}

\begin{corollary}\label{We}
Assume that $m\geq 1$ and $t\geq r$.  Letting $D= (\d_r, \d_t)$, one then has
\begin{equation}\label{We2}
|D^\b W_{jm}(\la,r,t)| \leq C(m)\cdot \left( \frac{r}{t+r} \right)^{1/2-|\b|} \cdot \frac{\la^{m-j+1/2}}{(t-r)^{m-j
+|\b|}} \cdot \frac{1}{\sqrt{t-r-\la}}
\end{equation}
for each $0\leq \la \leq t-r$ and each multi-index $\b$ with $|\b|\leq j\leq m$.
\end{corollary}

\begin{proof}
Let us first differentiate \eqref{W} to get
\begin{equation*}
D^\b W_{jm}(\la,r,t) = \int_{-1}^1 D^\b (\sigma-z)^{j-1/2} \cdot \frac{T_m(\sigma)}{\sqrt{1- \sigma^2}} \:\:d\sigma.
\end{equation*}
Using the chain rule repeatedly, one may write the leftmost factor in the integrand as
\begin{equation*}
D^\b (\sigma-z)^{j-1/2} = \sum_{s=1}^{|\b|} \:\sum_{\b_1+\ldots+\b_s= \b} C(\b_1,\ldots,\b_s,s,j) \cdot (\sigma -
z)^{j-s-1/2} \cdot \prod_{q=1}^s (D^{\b_q} z).
\end{equation*}
Recalling our definition \eqref{W}, one may then combine the last two equations to write
\begin{equation*}
D^\b W_{jm}(\la,r,t) = \sum_{s=1}^{|\b|} \:\sum_{\b_1+\ldots+\b_s= \b} C(\b_1,\ldots,\b_s,s,j) \cdot W_{j-s,m}(\la,r,t)
\cdot \prod_{q=1}^s (D^{\b_q} z).
\end{equation*}
This trivially gives
\begin{equation}\label{We3}
|D^\b W_{jm}(\la,r,t)| \leq C\sum_{s=1}^{|\b|} \:\sum_{\b_1+\ldots+\b_s= \b} |W_{j-s,m}(\la,r,t)| \cdot \prod_{q=1}^s
|D^{\b_q} z|
\end{equation}
and we shall now resort to our auxiliary estimate \eqref{dz3}.  It is easy to check that
\begin{equation*}
r^2 + |t^2-\la^2| \geq r^2 + r(t+\la) \geq rt, \quad\quad 0\leq \la\leq t-r
\end{equation*}
because $t-\la\geq r$ for such $\la$.  In particular, our auxiliary estimate \eqref{dz3} gives
\begin{equation*}
|D^{\b_q} z| \leq C(\b_q) \cdot \frac{r^2+t^2-\la^2}{\la r^{|\b_q|+1}}\leq \frac{C(\b_q)}{r^{|\b_q|}} \cdot \left(
\frac{t+r+\la}{\la} \cdot \frac{t+r-\la}{r} \right).
\end{equation*}
Employing this fact in \eqref{We3}, we thus arrive at
\begin{equation}\label{We4}
|D^\b W_{jm}(\la,r,t)| \leq Cr^{-|\b|} \cdot \sum_{s=1}^{|\b|} |W_{j-s,m}(\la,r,t)|  \cdot \left( \frac{t+r+\la}{\la}
\cdot \frac{t+r-\la}{r} \right)^s.
\end{equation}
Here, let us note that $0\leq |\b|-s \leq j-s \leq j\leq m$ since $|\b|\leq j\leq m$ by assumption.

\Case{1} When $0\leq \la \leq (t-r)/2$, it suffices to show that
\begin{equation}\label{g1}
|D^\b W_{jm}(\la,r,t)| \leq C(m)\cdot \left( \frac{r}{t+r} \right)^{1/2-|\b|} \cdot \frac{\la^{m-j+ 1/2}}{(t-r)^{m-j
+|\b|+1/2}}
\end{equation}
because $t-r-\la$ is equivalent to $t-r$.  For the exact same reason, we have
\begin{equation*}
|W_{j-s,m}(\la,r,t)| \leq C(m) \cdot \left( \frac{r}{t+r} \cdot \frac{\la}{t-r} \right)^{m-j+s+1/2}
\end{equation*}
by Lemma \ref{Wet} with $i=j-s$ and $i_0=m$.  Inserting this fact in \eqref{We4}, we then get
\begin{equation*}
|D^\b W_{jm}(\la,r,t)| \leq Cr^{-|\b|} \cdot \sum_{s=1}^{|\b|} \left( \frac{r}{t+r} \cdot \frac{\la}{t-r} \right)^{m-j
+1/2} \left( \frac{t+r+\la}{t+r} \cdot \frac{t+r-\la}{t-r} \right)^s.
\end{equation*}
Since $t+r\pm \la\leq 2(t+r)$ whenever $\la\leq t-r$, this also implies
\begin{align*}
|D^\b W_{jm}(\la,r,t)|
&\leq Cr^{-|\b|}\cdot \left( \frac{r}{t+r} \cdot \frac{\la}{t-r} \right)^{m-j+1/2} \left( \frac{t+r}{t-r} \right)^{|\b|}\\
&= C\left( \frac{r}{t+r}\right)^{m-j+1/2-|\b|} \cdot \frac{\la^{m-j+1/2}}{(t-r)^{m-j+|\b|+1/2}}\:.
\end{align*}
In view of our assumption that $j\leq m$, we may now deduce the desired estimate \eqref{g1}.

\Case{2} When $(t-r)/2 \leq \la \leq t-r$ and $r\leq t\leq 2r$, we need only show that
\begin{equation}\label{g2}
|D^\b W_{jm}(\la,r,t)| \leq \frac{C\la^{1/2-|\b|}}{\sqrt{t-r-\la}}\:.
\end{equation}
Since $t+r\pm \la \leq 2(t+r)\leq 6r$ under the present assumptions, we easily get
\begin{equation*}
|D^\b W_{jm}(\la,r,t)| \leq Cr^{-|\b|} \cdot \sum_{s=1}^{|\b|} |W_{j-s,m}(\la,r,t)| \cdot \left( \frac{r}{\la} \right)^s
\end{equation*}
by means of \eqref{We4}.  Moreover, $\la\leq t-r\leq r$ for this case, so we find that
\begin{equation*}
|D^\b W_{jm}(\la,r,t)| \leq C\la^{-|\b|} \cdot \sum_{s=1}^{|\b|} |W_{j-s,m}(\la,r,t)|.
\end{equation*}
Applying Lemma \ref{Wet} with $i=j-s=i_0$, we also have
\begin{equation*}
|W_{j-s,m}(\la,r,t)|\leq C(m)\cdot\left( \frac{r}{t+r} \cdot \frac{\la}{t-r-\la} \right)^{1/2} \leq
\frac{C\la^{1/2}}{\sqrt{t-r-\la}}
\end{equation*}
so we may simply combine the last two equations to deduce the desired estimate \eqref{g2}.

\Case{3} When $(t-r)/2 \leq \la \leq t-r$ and $t\geq 2r$, it suffices to show that
\begin{equation}\label{g3}
|D^\b W_{jm}(\la,r,t)| \leq \frac{Cr^{1/2-|\b|}}{\sqrt{t-r-\la}}
\end{equation}
because $\la$ is equivalent to $t\pm r$.  For the exact same reason, \eqref{We4} reduces to
\begin{equation*}
|D^\b W_{jm}(\la,r,t)| \leq Cr^{-|\b|} \cdot \sum_{s=1}^{|\b|} |W_{j-s,m}(\la,r,t)|  \cdot \left( \frac{t+r-\la}{r}
\right)^s,
\end{equation*}
and then Lemma \ref{Wet} applies to give
\begin{equation*}
|D^\b W_{jm}(\la,r,t)| \leq Cr^{-|\b|} \cdot \sum_{s=1}^{|\b|} \left( \frac{r}{t-r-\la} \right)^{i_0-j+s+1/2} \left(
\frac{t+r-\la}{r} \right)^s
\end{equation*}
for any integer $i_0$ with $j-s\leq i_0\leq m$.

\Subcase{3a} If it happens that $t-3r/2 \leq \la$, the admissible choice $i_0=j-s$ yields
\begin{equation*}
|D^\b W_{jm}(\la,r,t)| \leq Cr^{-|\b|} \cdot \sum_{s=1}^{|\b|} \left( \frac{r}{t-r-\la} \right)^{1/2} \left(
\frac{t+r-\la}{r} \right)^s \leq \frac{Cr^{1/2-|\b|}}{\sqrt{t-r-\la}}
\end{equation*}
because $t+r-\la\leq 5r/2$ whenever $t-3r/2 \leq \la$.

\Subcase{3b} If it happens that $\la \leq t-3r/2$, the admissible choice $i_0= j$ yields
\begin{equation*}
|D^\b W_{jm}(\la,r,t)| \leq Cr^{-|\b|} \cdot \sum_{s=1}^{|\b|} \left( \frac{r}{t-r-\la} \right)^{s+1/2} \left(
\frac{t+r-\la}{r} \right)^s \leq \frac{Cr^{1/2-|\b|}}{\sqrt{t-r-\la}}
\end{equation*}
because $t+r-\la \leq 5(t-r-\la)$ whenever $\la\leq t-3r/2$.

This establishes the desired \eqref{g3} and also concludes the proof of the corollary.
\end{proof}

We are finally in a position to prove an even-dimensional analogue of Lemma \ref{Lon}.

\begin{lemma}\label{Len}
Let $n\geq 4$ be an even integer and $a,m$ be as in \eqref{am}. Suppose $f_0\in \mathcal{C}^l(\R_+)$ and $f_1\in
\mathcal{C}^{l+1}(\R_+)$ for some integer $1\leq l\leq m$. Also, assume the singularity condition
\begin{equation}\label{sg}
\sum_{s=0}^{i+l-1} \la^s \,|f_i^{(s)}(\la)| = O \left( \la^{-2m-2+\de} \right) \quad \text{as $\la \rightarrow 0$}
\end{equation}
for $i=0,1$ and some fixed $\de>0$.  When $t\geq r$, one then has the analogue
\begin{align}\label{Le2}
\d_t^i [Lf_i](r,t) &= C_j \int_{t-r}^{t+r} [H_{ij}f_i](\la,r,t)\cdot r^{-m-a} U_{jm}(\la,r,t) \:d\la \notag \\
&\quad + C_j \int_0^{t-r} [H_{ij}f_i](\la,r,t)\cdot r^{-m-a} W_{jm}(\la,r,t) \:d\la
\end{align}
of \eqref{Lo2} for $i=0,1$ and all integers $i\leq j\leq l$. Here, the functions $U_{jm}$ and $W_{jm}$ are given by
\eqref{U} and \eqref{W}, respectively, while the linear operators $H_{ij}$ are the linear operators that we introduced in
Lemma \ref{Lon}.
\end{lemma}

\begin{proof}
For the case $i=0$ and $0\leq j\leq l$, we shall use induction on $j$ to verify the desired
\begin{align}\label{Le3}
[Lf_0](r,t) &= C_j \int_{t-r}^{t+r} [H_{0j}f_0](\la,r,t) \cdot r^{-m-a} U_{jm}(\la,r,t) \:d\la\notag\\
&\quad + C_j \int_0^{t-r} [H_{0j}f_0](\la,r,t) \cdot r^{-m-a} W_{jm}(\la,r,t) \:d\la.
\end{align}
When $j=0$, this follows by our various definitions \eqref{am}, \eqref{Le}, \eqref{H0j}, \eqref{U} and \eqref{W}.  Let us
now assume the last equation holds for some $j\geq 0$.  In view of the definition
\begin{equation*}
U_{jm}(\la,r,t) = \int_{z(\la,r,t)}^1 [\sigma-z(\la,r,t)]^{j-1/2} \cdot \frac{T_m(\sigma)}{\sqrt{1-\sigma^2}} \:\:d\sigma
\end{equation*}
we introduced in \eqref{U}, it is easy to see that $U_{jm}$ and $U_{j+1,m}$ are related by the formula
\begin{equation*}
U_{jm}(\la,r,t) = -\frac{1}{j+1/2} \cdot \frac{1}{\d_\la z(\la,r,t)} \cdot \d_\la U_{j+1,m}(\la,r,t).
\end{equation*}
Moreover, the same formula relates $W_{jm}$ and $W_{j+1,m}$ because of our definition \eqref{W}.  If we use this fact to
integrate by parts the integrals of \eqref{Le3}, we then get to replace the leftmost factor in each integrand by
\begin{equation*}
\frac{1}{j+1/2} \cdot \frac{\d}{\d\la} \left( \frac{1}{\d_\la z(\la,r,t)} \cdot [H_{0j} f_0](\la,r,t) \right) =
\frac{1}{j+1/2} \cdot [H_{0,j+1} f_0](\la,r,t)
\end{equation*}
in view of our definition \eqref{H0j}.  Since the integration by parts also replaces $U_{jm}$ by $U_{j+1,m}$ and $W_{jm}$
by $W_{j+1,m}$, we may thus complete the inductive proof of \eqref{Le3} by merely showing that no boundary terms arise in
the process. Now, the boundary term at $\la=t+r$ is zero because $U_{j+1,m}$ vanishes there by part (a) of Lemma
\ref{Ue}. At $\la= t-r$, we have two boundary terms with opposite signs.  In fact, the functions $U_{j+1,m}$ and
$W_{j+1,m}$ agree at that point by part (a) of Lemma \ref{Wet}, so these two boundary terms plainly cancel one another.
As for the boundary term at $\la= 0$, that one is given by
\begin{equation*}
\left[ -\frac{C_j}{j+1/2} \cdot \frac{[H_{0j}f_0](\la,r,t)}{\d_\la z(\la,r,t)} \cdot r^{-m-a} W_{j+1,m}(\la,r,t)
\right]_{\la=0}.
\end{equation*}
Using \eqref{dlz1}, \eqref{He} and \eqref{We2} together with our singularity condition \eqref{sg} on $f_0$, one can
easily check that this expression vanishes as well.  Thus, no boundary terms arise, indeed.

Finally, to establish \eqref{Le2} for the case $i=1$ and $1\leq j\leq l$, one proceeds as in the proof of Lemma
\ref{Lon}. Since only minor modifications are needed here, we shall omit the details.
\end{proof}

\section{The Riemann Operator in High Dimensions}\label{Rm}
Using the results of the previous section, we shall now study the Riemann operator, which is defined by either \eqref{Lo}
or \eqref{Le}, according to the parity of $n$.  To estimate its derivatives, we divide our analysis into two cases that
we treat separately in our next two propositions. First, we focus on the \textit{interior} region $t\geq 2r$ and provide
an estimate that does not involve any boundary terms. For odd values of $n$, a similar estimate appears in \cite{Kb1}.
The main idea we use to handle this region coincides with the one in \cite{PST1}, where the wave equation with the
inverse-square potential is studied. However, the approach in \cite{PST1} is quite different from ours.

\begin{prop}\label{Lin}
Let $n\geq 4$ be an integer and define $a,m$ by \eqref{am}. Suppose $f_0\in \mathcal{C}^l(\R_+)$ and $f_1\in
\mathcal{C}^{l+1}(\R_+)$ for some integer $1\leq l\leq m$.  Also, assume the singularity condition \eqref{sg} for $i=0,1$
and some fixed $\de>0$.  When $D= (\d_r,\d_t)$ and $t\geq 2r$, the Riemann operator is then such that
\begin{align}\label{Lin1}
|D^\b \d_t^i [Lf_i](r,t)| &\leq C_1(n) \cdot r^{j-|\b|-m-a} \int_{t-r}^{t+r} \frac{\la^{m-j}}{(r-t+\la)^{1-a}} \cdot
\sum_{s=0}^{i+j}
\la^{s+1-i} \,|f_i^{(s)}(\la)| \:d\la \notag \\
&\quad + \frac{C_2(n)\cdot r^{j-|\b|-m}}{(t-r)^{j+m+a}} \: \int_0^{t-r} \frac{\la^{2m}}{(t-r-\la)^{1-a}} \cdot
\sum_{s=0}^{i+j} \la^{s+1-i} \,|f_i^{(s)} (\la)| \:d\la
\end{align}
for $i=0,1$ and each $\max(i,|\b|) \leq j\leq l$.  Moreover, one has $C_2(n)= 0$ when $n$ is odd.
\end{prop}

\begin{proof}
We divide our analysis into two cases, according to the parity of $n$.

\Case{1} For odd values of $n$, an application of Lemma \ref{Lon} allows us to write
\begin{equation*}
\d_t^i [Lf_i](r,t) = \frac{(-1)^j}{2} \int_{t-r}^{t+r} [H_{ij}f_i](\la,r,t)\cdot r^{-m-a} P_{jm}(z(\la,r,t)) \:d\la,
\quad\quad i\leq j\leq l.
\end{equation*}
Recall that $P_{jm}$ is a polynomial with a zero of order $j$ at each of $\pm 1$, while $z(t\pm r,r,t) = \pm 1$ by
\eqref{z}.  In particular, we may differentiate the last equation $|\b|\leq j$ times without creating any boundary terms
to obtain
\begin{equation}\label{aa1}
|D^\b \d_t^i [Lf_i]| \leq C\int_{t-r}^{t+r} \sum_{\b_1+\b_2+\b_3=\b} \: |D^{\b_1}[H_{ij}f_i]| \cdot r^{-m-a-|\b_2|} \cdot
|D^{\b_3} P_{jm} (z(\la,r,t))| \:d\la.
\end{equation}
When it comes to the operators $H_{ij}$, the general estimate
\begin{equation}\label{ge}
\left| D^{\b_1} [H_{ij}f_i] \right| \leq \frac{C\la^j \,r^{j-|\b_1|} \,t^{|\b_1|-j}}{(t-r)^{i+j+|\b_1|}} \cdot
\sum_{s=0}^{i+j} \la^{m+a+s} \,|f_i^{(s)}(\la)|, \quad\quad 0\leq \la\leq t+r
\end{equation}
is provided by Lemma \ref{Lon}.  Under our assumption that $t\geq 2r$, it actually implies
\begin{equation}\label{ge2}
\left| D^{\b_1} [H_{ij}f_i] \right| \leq C\la^{-i-j} \,r^{j-|\b_1|} \cdot \sum_{s=0}^{i+j} \la^{m+a+s}
\,|f_i^{(s)}(\la)|, \quad\quad t-r\leq \la\leq t+r
\end{equation}
because $\la$, $t\pm r$ and $t$ are all equivalent here.  For the exact same reason, \eqref{dz5} gives
\begin{equation*}
|D^\c z(\la,r,t)| \leq C(\c) \cdot \frac{t+r}{\la r^{|\c|}} \leq \frac{C}{r^{|\c|}} \:,\quad\quad t-r\leq \la\leq t+r
\end{equation*}
for any multi-index $\c$.  Moreover, one has $|z(\la,r,t)|\leq 1$ for such $\la$, so this easily leads to
\begin{equation*}
|D^{\b_3} P_{jm} (z(\la,r,t))| \leq Cr^{-|\b_3|}.
\end{equation*}
Using the last equation and \eqref{ge2}, we may then estimate \eqref{aa1} as
\begin{equation*}
|D^\b \d_t^i [Lf_i](r,t)| \leq Cr^{j-|\b|-m-a} \int_{t-r}^{t+r}\la^{m-j} \cdot \sum_{s=0}^{i+j} \la^{s+a-i}
\,|f_i^{(s)}(\la)| \:d\la.
\end{equation*}
Since $a=1$ when $n$ is odd, this is precisely the desired estimate \eqref{Lin1} with $C_2(n)=0$.

\Case{2} For even values of $n$, an application of Lemma \ref{Len} allows us to write
\begin{align*}
\d_t^i [Lf_i](r,t) &= C_j \int_{t-r}^{t+r} [H_{ij}f_i](\la,r,t)\cdot r^{-m-a} U_{jm}(\la,r,t) \:d\la \notag \\
&\quad +C_j \int_0^{t-r} [H_{ij}f_i](\la,r,t)\cdot r^{-m-a} W_{jm}(\la,r,t) \:d\la, \quad\quad i\leq j\leq l.
\end{align*}
By Lemma \ref{Ue}, the function $U_{jm}$ vanishes with order $j$ at the endpoint $\la= t+r$.  Also, its derivatives up to
order $j$ agree with those of $W_{jm}$ at the endpoint $\la= t-r$ by Lemma \ref{Wet}. Thus, no boundary terms arise upon
$|\b|\leq j$ differentiations of the last equation, so we get
\begin{align}\label{aa2}
|D^\b \d_t^i [Lf_i]| &\leq C\int_{t-r}^{t+r} \sum_{\b_1+\b_2+\b_3= \b} |D^{\b_1} [H_{ij} f_i]| \cdot r^{-m-a-|\b_2|}
\cdot |D^{\b_3}
U_{jm}(\la,r,t)| \:d\la \notag \\
&\quad +C\int_0^{t-r} \sum_{\b_1+\b_2+\b_3= \b} |D^{\b_1} [H_{ij} f_i]| \cdot r^{-m-a-|\b_2|} \cdot |D^{\b_3}
W_{jm}(\la,r,t)| \:d\la \notag \\
&\equiv A_1 + A_2.
\end{align}
Within the first integral, $\la$ is equivalent to $t\pm r$ as before, so an estimate of the form
\begin{equation*}
|D^{\b_3} U_{jm}(\la,r,t)| \leq \frac{Cr^{1/2-|\b_3|}}{\sqrt{r-t+\la}} \:,\quad\quad t-r\leq \la\leq t+r
\end{equation*}
holds by part (c) of Lemma \ref{Ue}.  Using this fact along with \eqref{ge2}, we then find that
\begin{equation*}
A_1 \leq Cr^{j-|\b|-m-a+1/2} \int_{t-r}^{t+r} \frac{\la^{m-j}}{\sqrt{r-t+\la}} \cdot \sum_{s=0}^{i+j} \la^{s+a-i}
\,|f_i^{(s)}(\la)| \:d\la.
\end{equation*}
Moreover, we have $r\leq t-r\leq \la$ within the region of integration, hence also
\begin{equation*}
A_1 \leq Cr^{j-|\b|-m-a} \int_{t-r}^{t+r} \frac{\la^{m-j}}{\sqrt{r-t+\la}} \cdot \sum_{s=0}^{i+j} \la^{s+a+1/2-i}
\,|f_i^{(s)}(\la)| \:d\la.
\end{equation*}
Since $a=1/2$ when $n$ is even, the desired estimate \eqref{Lin1} is thus satisfied by $A_1$.

To treat the second integral $A_2$, we recall that $t$ is equivalent to $t\pm r$ because $t\geq 2r$ by assumption. An
immediate consequence of \eqref{ge} is then
\begin{equation*}
|D^{\b_1} [H_{ij}f_i]| \leq \frac{C\la^j \,r^{j-|\b_1|}}{(t-r)^{i+2j}} \cdot \sum_{s=0}^{i+j} \la^{m+a+s}
\,|f_i^{(s)}(\la)|, \quad\quad 0\leq \la\leq t-r
\end{equation*}
and an immediate consequence of Corollary \ref{We} is
\begin{equation*}
|D^{\b_3} W_{jm}(\la,r,t)| \leq \frac{Cr^{1/2-|\b_3|}\cdot \la^{m-j+1/2}}{(t-r)^{m-j+1/2}} \cdot \frac{1}{\sqrt{t-r-\la}}
\:,\quad\quad 0\leq \la\leq t-r.
\end{equation*}
Employing these two facts in \eqref{aa2}, we now get
\begin{equation*}
A_2 \leq \frac{Cr^{j-|\b|-m-a+1/2}}{(t-r)^{i+j+m+1/2}} \cdot \int_0^{t-r} \frac{\la^{2m}}{\sqrt{t-r-\la}} \cdot
\sum_{s=0}^{i+j} \la^{s+a+1/2} \,|f_i^{(s)}(\la)| \:d\la.
\end{equation*}
Moreover, $\la\leq t-r$ within the region of integration, so this also implies
\begin{equation*}
A_2 \leq \frac{Cr^{j-|\b|-m-a+1/2}}{(t-r)^{j+m+1/2}} \cdot \int_0^{t-r} \frac{\la^{2m}}{\sqrt{t-r-\la}} \cdot
\sum_{s=0}^{i+j} \la^{s+a+1/2-i} \,|f_i^{(s)}(\la)| \:d\la.
\end{equation*}
Since $a=1/2$ when $n$ is even, the desired estimate \eqref{Lin1} is thus satisfied by $A_2$ as well.
\end{proof}

In our next proposition, we provide a similar analysis for the \textit{exterior} region $t\leq 2r$.  The most notable
difference from our previous estimate is the presence of some boundary terms. A more subtle difference is that the
regularity of the initial data is not so crucial anymore.  Namely, our sharpest conclusions occur for the smallest
possible value of $j$ in the following

\begin{prop}\label{Lex}
When $t\leq 2r$, the assumptions of the previous proposition imply
\begin{align}\label{Lex1}
|D^\b \d_t^i [Lf_i](r,t)| &\leq C_1'(n) \cdot r^{-m-a} \int_{|t-r|}^{t+r} \frac{\la^{m-|\b|}}{(r-t+\la)^{1-a}} \cdot
\sum_{s=0}^{i+j} \la^{s+1-i}
\,|f_i^{(s)}(\la)| \:d\la \notag \\
&\quad + \frac{C_2'(n)\cdot r^{-m-a}}{(t-r)^{m+|\b|}} \:\int_0^{\max(t-r,0)} \frac{\la^{2m}}{(t-r-\la)^{1-a}} \cdot
\sum_{s=0}^{i+j} \la^{s+1-i} \,|f_i^{(s)}(\la)| \:d\la \notag\\
&\quad + C_1'(n)\cdot r^{-m-a} \: \sum_{s=0}^{i+j-1} \left[ \la^{m+a-|\b|+s+1-i} \,|f_i^{(s)}(\la)| \right]_{\la= |t\pm
r|}
\end{align}
for $i=0,1$ and each $\max(i,|\b|) \leq j\leq l$.  Moreover, one has $C_2'(n)=0$ when $n$ is odd.
\end{prop}

\begin{proof}
We divide our analysis into three cases.

\Case{1} For odd values of $n$, the Riemann operator \eqref{Lo} is of the form
\begin{equation}\label{Lg1}
[Lf_i](r,t) = \frac{1}{2} \int_{|t-r|}^{t+r} \la^{m+a} f_i(\la) \cdot r^{-m-a} P_m(z(\la,r,t)) \:d\la.
\end{equation}
To obtain the desired estimate in this case, we shall not have to distinguish between radial and time derivatives, so it
is convenient to introduce a multi-index $\c$ of order $|\c| = |\b|+i$.  Once we differentiate the last equation, we then
get an identity of the form
\begin{align}\label{bb}
D^\c [Lf_i](r,t)
&= \frac{1}{2} \int_{|t-r|}^{t+r} \la^{m+a} f_i(\la) \cdot D^\c \Bigl( r^{-m-a} P_m(z(\la,r,t)) \Bigr) \:d\la \notag \\
& \quad + \sum_{|\c_1|= |\c|-1} C_\pm(\c_1) \cdot D^{\c_1} \Bigl[ \la^{m+a} f_i(\la) \cdot
r^{-m-a} P_m(z(\la,r,t)) \Bigr]_{\la= |t\pm r|} \notag \\
&\equiv B_1 + B_2.
\end{align}
Since $|t-r|\leq \la\leq t+r$ within these terms, our assumption $t\leq 2r$ makes Corollary \ref{dzb} applicable. Letting
$D_*= (\d_\la, \d_r, \d_t)$, we thus have the estimate
\begin{equation*}
|D_*^\al z(\la,r,t)| \leq C(\al)\cdot \la^{-|\al|}, \quad\quad |t-r|\leq \la\leq t+r
\end{equation*}
for any multi-index $\al$.  Besides, $|z(\la,r,t)|\leq 1$ for such $\la$, so we easily get
\begin{equation}\label{Pmd}
|D_*^\al P_m(z(\la,r,t))| \leq C(\al,m) \cdot \la^{-|\al|}, \quad\quad |t-r|\leq \la\leq t+r
\end{equation}
by repeated applications of the chain rule.  This is actually the only fact we need to control the right hand side of
\eqref{bb}. Due to our assumption that $t\leq 2r$, it also implies
\begin{equation}\label{Pm2}
\left| D_*^\al \left( r^{-m-a} P_m(z(\la,r,t)) \right) \right| \leq Cr^{-m-a} \cdot \la^{-|\al|}, \quad\quad |t-r|\leq
\la\leq t+r
\end{equation}
for any multi-index $\al$, as $\la\leq 3r$ here.  Applying the last inequality to the integral term $B_1$, let us now
recall that $|\c|=|\b|+i$ to find
\begin{align*}
|B_1| \leq Cr^{-m-a} \int_{|t-r|}^{t+r} \la^{m+a-|\c|} \cdot |f_i(\la)| \:d\la = Cr^{-m-a} \int_{|t-r|}^{t+r}
\la^{m-|\b|} \cdot \la^{a-i} \,|f_i(\la)| \:d\la.
\end{align*}
Since $a=1$ when $n$ is odd, the desired estimate \eqref{Lex1} thus follows for $B_1$.  Applying \eqref{Pm2} to the
boundary terms $B_2$, we similarly get
\begin{align*}
|B_2| &\leq \sum_{|\c_1|= |\b|+i-1} \:\sum_{s_1+s_2= |\c_1|} Cr^{-m-a} \cdot \Bigl[ \la^{m+a-s_1} \,|f_i^{(s_2)}(\la)|
\Bigr]_{\la= |t\pm r|} \\
&= Cr^{-m-a} \: \sum_{s_2=0}^{i+|\b|-1} \Bigl[ \la^{m+a-|\b|+s_2+1-i} \,|f_i^{(s_2)}(\la)| \Bigr]_{\la= |t\pm r|}.
\end{align*}
Moreover, $|\b|\leq j$ by assumption, so the desired estimate \eqref{Lex1} follows for $B_2$ as well.

\Case{2} When $n$ is even and $t\leq r$, the Riemann operator \eqref{Le} takes the form
\begin{equation}\label{Lg2}
[Lf_i](r,t) = C_0' \int_{|t-r|}^{t+r} \la^{m+a} f_i(\la) \cdot r^{-m-a} U_{0m}(\la,r,t) \:d\la
\end{equation}
with $U_{0m}$ given by \eqref{U}.  This closely resembles the equation \eqref{Lg1} we used in the previous case, although
the factor $P_m(z)$ has been replaced by $U_{0m}$ and the value of $a$ has changed due to the change in the parity of
$n$. Differentiating directly as before, we now get
\begin{equation*}
D^\c [Lf_i](r,t) = B_1' + B_2',
\end{equation*}
where each $B_k'$ denotes the corresponding $B_k$ of \eqref{bb} with $U_{0m}$ instead of $P_m(z)$.

Since $t\leq 2r$ by assumption, part (b) of Lemma \ref{Ue} applies to give
\begin{equation}\label{ig}
|D_*^\al U_{0m}(\la,r,t)| \leq C\la^{-|\al|} \cdot \left( \frac{\la}{r-t+\la} \right)^{1/2+|\al|}, \quad\quad |t-r|\leq
\la\leq t+r
\end{equation}
for any multi-index $\al$. This provides the first fact that we are going to need; the second one is the inequality
$\la\leq r-t+\la$, which trivially holds for the case $t\leq r$ under consideration.  Combining these two facts, we see
that \eqref{Pmd} remains valid when $P_m(z)$ is replaced by $U_{0m}$.  In particular, the argument of Case 1 applies
verbatim, except for the part where the exact value of $a$ was invoked. Said differently, we need only worry about the
integral term
\begin{equation}\label{B1'}
B_1' = C_0' \int_{|t-r|}^{t+r} \la^{m+a} f_i(\la) \cdot D^\c \Bigl( r^{-m-a} U_{0m}(\la,r,t) \Bigr) \:d\la.
\end{equation}
If we combine \eqref{ig} with our inequality $\la\leq r-t+\la$ once again, we get
\begin{equation*}
\left| D_*^\al U_{0m}(\la,r,t) \right| \leq \frac{C\la^{1/2-|\al|}}{\sqrt{r-t+\la}}\:, \quad\quad |t-r|\leq \la\leq t+r
\end{equation*}
for any multi-index $\al$.  As long as $\la\leq t+r\leq 3r$, we then get
\begin{equation*}
\left| D_*^\al \left( r^{-m-a} U_{0m}(\la,r,t) \right) \right| \leq \frac{Cr^{-m-a} \cdot
\la^{1/2-|\al|}}{\sqrt{r-t+\la}} \:,\quad\quad |t-r|\leq \la\leq t+r
\end{equation*}
for any multi-index $\al$. Applying this fact to the integral term \eqref{B1'}, we now find
\begin{align*}
|B_1'| \leq Cr^{-m-a} \int_{|t-r|}^{t+r} \frac{\la^{1/2-|\c|}}{\sqrt{r-t+\la}} \cdot \la^{m+a} \,|f_i(\la)| \:d\la.
\end{align*}
Since $a=1/2$ when $n$ is even and $|\c|= |\b|+i$ by above, we thus obtain
\begin{align*}
|B_1'| \leq Cr^{-m-a} \int_{|t-r|}^{t+r} \frac{\la^{m-|\b|}}{(r-t+\la)^{1-a}} \cdot \la^{1-i} \,|f_i(\la)| \:d\la.
\end{align*}
In particular, we obtain the desired estimate \eqref{Lex1} for this case as well.

\Case{3} Suppose now that $n$ is even and $r<t\leq 2r$.  Although our estimate \eqref{ig} behaves nicely at the boundary
terms $\la= r\pm t$ we had before, this is not the case for the boundary term $\la= t-r$ which may now emerge.  One may
overcome this difficulty using our approach in the previous proposition to avoid the boundary terms altogether. However,
that approach does not lead to the desired estimate, so we shall need to resort to a combination of the two.

Proceeding to the details, let us fix some $\zeta \in \mathcal{C}^\infty(\R)$ such that
\begin{equation}\label{zi1}
\zeta(x) = \left\{ \begin{array}{clll} 1 &&\text{if\, $x\leq 2$}
\\ 0 &&\text{if\, $x\geq 5/2$} \end{array} \right\}
\end{equation}
and introduce the functions
\begin{equation}\label{zi2}
\zeta_1(\la) = 1- \zeta \left( \frac{\la}{t-r} \right), \quad \zeta_2(\la) = \zeta \left( \frac{\la}{t-r} \right).
\end{equation}
Then the linearity of the Riemann operator allows us to write
\begin{equation*}
[Lf_i](r,t) = [L(\zeta_1 f_i)](r,t) + [L(\zeta_2 f_i)](r,t)
\end{equation*}
as a sum of two terms to be treated separately.

In order to treat $[L(\zeta_1 f_i)]$, we proceed as in Case 2. Since $\zeta_1$ vanishes when $\la\leq 2(t-r)$ by
definition, the Riemann operator \eqref{Le} takes the form
\begin{equation*}
[L(\zeta_1 f_i)](r,t) = C_0' \int_{2(t-r)}^{t+r} \la^{m+a} \zeta_1(\la)f_i(\la) \cdot r^{-m-a} U_{0m}(\la,r,t) \:d\la.
\end{equation*}
This closely resembles the equation \eqref{Lg2} we had in Case 2, although the extra factor $\zeta_1(\la)$ is now present
in the integrand.  When it comes to the first fact \eqref{ig} we used in Case 2, it is provided by Lemma \ref{Ue}
whenever $t\leq 2r$, so it is still valid.  As for the inequality we used, its analogue $\la\leq 2(r-t+\la)$ now holds
since $\la\geq 2(t-r)$ within the region of integration.  The main question is then whether the extra factor
$\zeta_1(\la)$ will bring any changes to our previous approach. Although this function is bounded by our definition
\eqref{zi2}, its derivatives involve powers of $t-r$ we did not have before. Nevertheless, the only place where such
derivatives may occur is at the endpoint $\la= t+r$, and one has $t+r\geq 3(t-r)$ whenever $t\leq 2r$.  Once we now
recall our definition \eqref{zi1}-\eqref{zi2}, we see that all derivatives of $\zeta_1$ vanish at $t+r$.  In particular,
one does not have to deal with such derivatives at all.

Next, we focus on $[L(\zeta_2 f_i)]$.  For $\la$ sufficiently close to zero, $\zeta_2(\la)\equiv 1$ by
\eqref{zi1}-\eqref{zi2}, so the functions $\zeta_2 f_i$ and $f_i$ satisfy the same singularity condition as $\la\to 0$.
This means that we may employ Lemma \ref{Len} and argue as in the previous proposition to find
\begin{align}\label{bb1}
|D^\b \d_t^i [L(\zeta_2f_i)]| &\leq C\int_{t-r}^{t+r} \sum_{\b_1+\b_2+\b_3= \b} |D^{\b_1} [H_{ij} (\zeta_2f_i)]| \cdot
r^{-m-a-|\b_2|} \cdot |D^{\b_3}
U_{jm}| \:d\la \notag \\
&\quad + C\int_0^{t-r} \sum_{\b_1+\b_2+\b_3= \b} |D^{\b_1} [H_{ij} (\zeta_2f_i)]| \cdot r^{-m-a-|\b_2|} \cdot |D^{\b_3}
W_{jm}| \:d\la \notag \\
&\equiv B_1'' + B_2''
\end{align}
for $i=0,1$ and each $\max(i,|\b|)\leq j\leq l$.  Namely, the condition $i\leq j\leq l$ allows us to invoke the identity
of Lemma \ref{Len} and the condition $|\b|\leq j$ to differentiate that identity $|\b|$ times without introducing any
boundary terms.

Let us note that $3(t-r)\leq t+r$ whenever $t\leq 2r$ and temporarily assume the estimate
\begin{equation}\label{ass}
|D^{\b_1} [H_{ij}(\zeta_2f_i)]| \leq \frac{C\la^j}{(t-r)^{i+j+|\b_1|}} \cdot \sum_{s_2=0}^{i+j} \la^{m+a+s_2}
\,|f_i^{(s_2)}(\la)|, \quad\quad 0\leq \la\leq 3(t-r)
\end{equation}
whose proof is given below.  For the remaining values $\la\geq 3(t-r)$, the left hand side is zero because $\zeta_2$
vanishes for such $\la$ by \eqref{zi1}-\eqref{zi2}. Employing \eqref{ass}, one thus arrives at
\begin{equation*}
B_1'' \leq Cr^{-m-a} \int_{t-r}^{3(t-r)} \sum_{\b_1+\b_2+\b_3= \b} r^{-|\b_2|} \cdot |D^{\b_3} U_{jm}| \cdot
\la^{-i-|\b_1|} \sum_{s_2=0}^{i+j} \la^{m+a+s_2} \,|f_i^{(s_2)}(\la)|\:d\la.
\end{equation*}
Since $r$ is equivalent to $t+r$ whenever $t\leq 2r$, part (c) of Lemma \ref{Ue} ensures that
\begin{equation*}
|D^{\b_3} U_{jm}(\la,r,t)| \leq \frac{C\la^{1/2-|\b_3|}}{\sqrt{r-t+\la}} \:,\quad\quad t-r\leq \la\leq 3(t-r)
\end{equation*}
because $3(t-r)\leq t+r$ by above.  Once we now combine the last two equations, we get
\begin{equation*}
B_1'' \leq Cr^{-m-a} \int_{t-r}^{3(t-r)} \sum_{\b_1+\b_2+\b_3= \b} \frac{r^{-|\b_2|} \,\la^{m-|\b_1|-|\b_3|}}{\sqrt{r-t
+\la}} \cdot \sum_{s_2=0}^{i+j} \la^{s_2+a+1/2-i} \,|f_i^{(s_2)}(\la)|\:d\la.
\end{equation*}
Moreover, $r\geq t-r\geq \la/3$ within the region of integration, so this also gives
\begin{equation*}
B_1'' \leq Cr^{-m-a} \int_{t-r}^{t+r} \frac{\la^{m-|\b|}}{\sqrt{r-t+\la}} \cdot \sum_{s_2=0}^{i+j} \la^{s_2+a+1/2-i}
\,|f_i^{(s_2)}(\la)|\:d\la.
\end{equation*}
Since $a=1/2$ when $n$ is even, the desired estimate \eqref{Lex1} is thus satisfied by $B_1''$.

Next, we focus on $B_2''$.  Noting that $r\leq t+r\leq 3r$ by assumption, we get
\begin{equation*}
|D^{\b_3} W_{jm}(\la,r,t)| \leq \frac{C\la^{m-j+1/2}}{(t-r)^{m-j+|\b_3|}}\cdot \frac{1}{\sqrt{t-r-\la}} \:,\quad\quad
0\leq \la\leq t-r
\end{equation*}
by means of Corollary \ref{We}.  Together with \eqref{ass}, this allows us to estimate $B_2''$ as
\begin{equation*}
B_2'' \leq \frac{Cr^{-m-a}}{(t-r)^{i+m}} \int_0^{t-r} \sum_{\b_1+\b_2+\b_3= \b} \frac{r^{-|\b_2|}}{(t-r)^{|\b_1|+
|\b_3|}} \cdot \frac{\la^{m+1/2}}{\sqrt{t-r-\la}} \cdot \sum_{s_2=0}^{i+j} \la^{m+a+s_2} \,|f_i^{(s_2)}(\la)| \:d\la.
\end{equation*}
Moreover, $r^{-|\b_2|}(t-r)^{-i} \leq (t-r)^{-|\b_2|}\la^{-i}$ whenever $\la\leq t-r\leq r$, so we find
\begin{equation*}
B_2'' \leq \frac{Cr^{-m-a}}{(t-r)^{m+|\b|}} \int_0^{t-r} \frac{\la^{2m}}{\sqrt{t-r-\la}} \cdot \sum_{s_2=0}^{i+j}
\la^{s_2+a+1/2-i} \,|f_i^{(s_2)}(\la)| \:d\la.
\end{equation*}
Since $a=1/2$ when $n$ is even, the desired estimate \eqref{Lex1} is thus satisfied by $B_2''$ as well.

To finish the proof, it remains to establish \eqref{ass}. According to Lemma \ref{Lon}, we do have
\begin{equation*}
|D^{\b_1} [H_{ij}(\zeta_2f_i)]| \leq \frac{C\la^j \,r^{j-|\b_1|} \,t^{|\b_1|-j}}{(t-r)^{i+j+|\b_1|}} \cdot
\sum_{s=0}^{i+j} \la^{m+a+s} \cdot \left| \frac{d^s}{d\la^s} \,(\zeta_2(\la)f_i(\la)) \right|
\end{equation*}
for each $0\leq \la \leq 3(t-r)$ because $3(t-r)\leq t+r$ by above.  For the case $r< t\leq 2r$ we are presently
considering, our remaining assertion \eqref{ass} will then follow once we show that
\begin{equation}\label{sh}
\sum_{s=0}^{i+j} \la^s \cdot \left| \frac{d^s}{d\la^s} \,(\zeta_2(\la)f_i(\la)) \right| \leq C \sum_{s_2=0}^{i+j}
\la^{s_2} \cdot |f_i^{(s_2)}(\la)|, \quad\quad 0\leq \la\leq 3(t-r).
\end{equation}
Recalling our definition \eqref{zi1}-\eqref{zi2}, we easily get this when $2(t-r)\leq \la\leq 3(t-r)$ since
\begin{equation*}
\sum_{s=0}^{i+j} \la^s \cdot \left| \frac{d^s}{d\la^s} \,(\zeta_2(\la)f_i(\la)) \right| \leq C \sum_{s=0}^{i+j}
\sum_{s_2=0}^s  \la^s (t-r)^{s_2-s} \cdot |f_i^{(s_2)}(\la)| \leq C \sum_{s_2=0}^{i+j} \la^{s_2} \cdot |f_i^{(s_2)}(\la)|
\end{equation*}
by Leibniz' rule.  When $0\leq \la< 2(t-r)$, on the other hand, our definition \eqref{zi1}-\eqref{zi2} is such that all
derivatives of $\zeta_2$ vanish at $\la$, whence
\begin{equation*}
\sum_{s=0}^{i+j} \la^s \cdot \left| \frac{d^s}{d\la^s} \left( \zeta_2(\la) f_i(\la) \right) \right| = \sum_{s=0}^{i+j}
\la^s \cdot |\zeta_2(\la)| \cdot |f_i^{(s)}(\la)| \leq C \sum_{s=0}^{i+j} \la^s \cdot |f_i^{(s)}(\la)|.
\end{equation*}
In any case whatsoever, the desired estimate \eqref{sh} follows and the proof is complete.
\end{proof}

\section{Estimates for the Free Solution}\label{fr}
Using the results of the previous section, we now study the solution $u_0$ of the homogeneous equation \eqref{he}.  To
prove the estimates that Theorem \ref{dec} asserts for its derivatives, we shall also need two elementary facts that we
list separately in our next two lemmas.  Their proofs follow those of similar results in \cite{As, Ts1}, for instance, so
we are going to omit them. As for the parameter $a>0$ we introduce below, we are merely interested in the special cases
$a=1$ and $a=1/2$ that we had before \eqref{am}.

\begin{lemma}\label{ts1}
Let $(r,t)\in \R_+^2$ be arbitrary.  Assuming that $a>0$, one has
\begin{equation*}
\int_{|t-r|}^{t+r} \frac{\br{\la}^b d\la}{(r-t+\la)^{1-a}} \leq \left\{
\begin{array}{llcc}
Cr^a \br{t+r}^b &&\text{if} & b > -a \\
Cr^a \br{t+r}^{-a} \left( 1+ \ln \dfrac{\br{t+r}}{\br{t-r}} \right) &&\text{if} & b = -a \\
Cr^a \br{t+r}^{-a} \br{t-r}^{a+b} &&\text{if} & b < -a \\
\end{array}\right\}
\end{equation*}
for some constant $C$ depending solely on $a$ and $b$.
\end{lemma}

\begin{lemma}\label{ts2}
Let $a>0$ be arbitrary.  Assuming that $b\geq 0$ and $t\geq r>0$, the integral
\begin{equation*}
\mathcal{B}(b,c) \equiv \int_0^{t-r} \frac{\la^b \br{\la}^c \:d\la}{(t-r-\la)^{1-a}}
\end{equation*}
satisfies an estimate of the form
\begin{equation*}
\mathcal{B}(b,c) \leq \left\{
\begin{array}{lllll}
C(t-r)^{a+b} \cdot \br{t-r}^c &&\text{if} &b+c> -1 \\
C(t-r)^{a+b} \cdot \br{t-r}^{-b-1}\cdot (1+ \ln \br{t-r}) &&\text{if} &b+c= -1\\  C(t-r)^{a+b} \cdot \br{t-r}^{-b-1}
&&\text{if}& b+c< -1
\end{array}\right\},
\end{equation*}
where the constant $C$ depends solely on $a,b$ and $c$.
\end{lemma}

\begin{proof_of}{Theorem \ref{dec}}
Let $u_0(r,t)$ be the solution of the homogeneous equation \eqref{he} provided by Lemma \ref{hs}.  Setting $f_0= \psi$
and $f_1= \phi$, we may then write
\begin{equation*}
D^\b u_0(r,t) = D^\b [L\psi](r,t) + D^\b \d_t [L\phi](r,t) = \sum_{i=0}^1 D^\b \d_t^i [Lf_i](r,t).
\end{equation*}
This is also the expression we wish to estimate when $|\b|\leq l$. Using the notation above, let us now write our
assumption \eqref{data} on the initial data simply as
\begin{equation}\label{da2}
\sum_{i=0}^1 \sum_{s=0}^{i+l} \la^{s+1-i} \,|f_i^{(s)}(\la)| \leq \e\la^{l-m} \br{\la}^{m-l-k}.
\end{equation}
Here, one may readily check that the singularity condition \eqref{sg} holds for $i=0,1$.

\Step{1} In the interior region $t\geq 2r$, we use Proposition \ref{Lin} with $j= l$ to obtain
\begin{align*}
|D^\b u_0(r,t)| &\leq C_1 r^{l-|\b|-m-a} \int_{t-r}^{t+r} \frac{\la^{m-l}}{(r-t+\la)^{1-a}} \cdot \sum_{i=0}^1
\sum_{s=0}^{i+l}
\la^{s+1-i} \,|f_i^{(s)}(\la)| \:d\la \\
&\quad + \frac{C_2 r^{l-|\b|-m}}{(t-r)^{l+m+a}} \: \int_0^{t-r} \frac{\la^{2m}}{(t-r-\la)^{1-a}} \cdot \sum_{i=0}^1
\sum_{s=0}^{i+l} \la^{s+1-i} \,|f_i^{(s)} (\la)| \:d\la,
\end{align*}
where $C_2 = 0$ for odd values of $n$.  Once we now combine this with \eqref{da2}, we get
\begin{align*}
|D^\b u_0| &\leq C_1\e r^{l-|\b|-m-a} \int_{t-r}^{t+r} \frac{\br{\la}^{m-l-k} d\la}{(r-t+\la)^{1-a}}
+\frac{C_2\e r^{l-|\b| -m}}{(t-r)^{l+m+a}} \: \int_0^{t-r} \frac{\la^{l+m}\br{\la}^{m-l-k} d\la}{(t-r-\la)^{1-a}} \\
&\equiv \mathbb{A}_1 + \mathbb{A}_2,
\end{align*}
where $\mathbb{A}_2 = 0$ for odd values of $n$.

In what follows, we only deal with the decay rates $0\leq k < 2(m+a)$ that correspond to parts (a) through (c), omitting
the similar approach that one needs for parts (d) and (e).  Since $t+r$ and $t-r$ are equivalent whenever $t\geq 2r$, an
inequality of the form
\begin{equation}\label{go1}
|D^\b u_0| \leq C\e r^{l-|\b|-m} \cdot \br{t+r}^{m-l-k}
\end{equation}
implies each of the three desired estimates.  One can easily obtain this upper bound for $\mathbb{A}_1$ since $\la$ is
equivalent to $t\pm r$ within the region of integration.  As for the integral term $\mathbb{A}_2$, we may assume that $n$
is even and thus focus on decay rates $0\leq k< 2m+1$. Write
\begin{equation*}
\mathbb{A}_2 = C_2\e r^{l-|\b|-m} \cdot (t-r)^{-l-m-a} \cdot \mathcal{B}(l+m, m-l-k)
\end{equation*}
using the notation of Lemma \ref{ts2}.  Since $2m-k>-1$ by above, the lemma ensures that
\begin{equation}\label{B3}
\mathcal{B}(l+m, m-l-k) \leq C(t-r)^{l+m+a} \cdot \br{t-r}^{m-l-k}.
\end{equation}
Once we now combine the last two equations, we may deduce the desired \eqref{go1}.

\Step{2} In the exterior region $t\leq 2r$, Proposition \ref{Lex} with $j= \max(i,|\b|)$ similarly gives
\begin{align*}
|D^\b u_0(r,t)| &\leq C_1'\e r^{-m-a} \int_{|t-r|}^{t+r} \frac{\la^{l-|\b|} \br{\la}^{m-l-k} d\la}{(r-t +\la)^{1-a}} +
C_1'\e r^{-m-a} \: |t\pm r|^{a+l-|\b|} \br{t\pm r}^{m-l-k} \notag \\
&\quad + \frac{C_2'\e r^{-m-a}}{(t-r)^{m+|\b|}}\int_0^{\max(t-r,0)} \frac{\la^{l+m} \br{\la}^{m-l-k}\:d\la}{(t-r
-\la)^{1-a}} \notag \\
&\equiv \mathbb{B}_1 + \mathbb{B}_2 + \mathbb{B}_3,
\end{align*}
where $\mathbb{B}_3=0$ for odd values of $n$.  Let us only consider the decay rates $0\leq k< 2(m+a)$, as before. In view
of \eqref{B3}, the integral term $\mathbb{B}_3$ is then such that
\begin{equation*}
\mathbb{B}_3 \leq C\e r^{-m-a} \cdot |t-r|^{a+l-|\b|} \cdot \br{t-r}^{m-l-k}.
\end{equation*}
Since the right hand side appears in the boundary terms $\mathbb{B}_2$, it thus suffices to treat
\begin{equation}\label{B12}
\mathbb{B}_1 + \mathbb{B}_2 = C_1'\e r^{-m-a} \left[ \int_{|t-r|}^{t+r} \frac{\la^{l-|\b|} \br{\la}^{m-l-k} d\la}{(r-t
+\la)^{1-a}} + |t\pm r|^{a+l-|\b|} \br{t\pm r}^{m-l-k} \right].
\end{equation}

\Case{1} If $t\leq 2r$ and $r\leq 1$, then $t+r$ is bounded and we need only show that
\begin{equation*}
|D^\b u_0(r,t)| \leq C\e r^{l-|\b|-m}.
\end{equation*}
In fact, $|t\pm r| \leq 3r$ for this case, so we easily get the desired estimate
\begin{equation*}
\mathbb{B}_1 + \mathbb{B}_2 \leq C\e r^{l-|\b|-m-a} \left[ \int_{|t-r|}^{t+r} (r-t +\la)^{a-1} \:d\la + r^a \right] \leq
C\e r^{l-|\b|-m}
\end{equation*}
because $|\b|\leq l$ and $a>0$ by assumption.

\Case{2} Suppose now that $t\leq 2r$ and $r\geq 1$.  Since $|\b|\leq l$ by assumption, \eqref{B12} leads to
\begin{align*}
\mathbb{B}_1 + \mathbb{B}_2 &\leq C\e r^{-m-a} \left[ \int_{|t-r|}^{t+r} \frac{\br{\la}^{m-k-|\b|}
\:d\la}{(r-t+\la)^{1-a}} + \br{t\pm r}^{m+a-k-|\b|} \right] \\
&\leq C\e r^{-m-a} \br{t-r}^{-|\b|} \left[ \int_{|t-r|}^{t+r} \frac{\br{\la}^{m-k} \:d\la}{(r-t+\la)^{1-a}} + \br{t\pm
r}^{m+a-k} \right]
\end{align*}
and we may use Lemma \ref{ts1} to handle the integral.

\Subcase{2a} When $0\leq k< m+a$, we find that
\begin{equation*}
\mathbb{B}_1 + \mathbb{B}_2 \leq C\e r^{-m-a} \br{t-r}^{-|\b|} \left[ r^a \br{t+r}^{m-k} + \br{t\pm r}^{m+a-k} \right].
\end{equation*}
Moreover, $m+a-k$ is positive, so we get
\begin{equation*}
\mathbb{B}_1 + \mathbb{B}_2 \leq C\e r^{-m-a} \br{t-r}^{-|\b|} \br{t+r}^{m+a-k}.
\end{equation*}
This does imply the desired estimate, as $r$ is equivalent to $\br{t+r}$ when $r\geq \max(t/2,1)$.

\Subcase{2b} When $k=m+a$, our previous approach applies verbatim, although an extra logarithmic factor is now included
in the estimate that Lemma \ref{ts1} provides.

\Subcase{2c} When $m+a< k< 2(m+a)$, Lemma \ref{ts1} yields
\begin{align*}
\mathbb{B}_1 + \mathbb{B}_2
&\leq C\e r^{-m-a} \br{t-r}^{-|\b|} \left[ r^a \br{t+r}^{-a} \br{t-r}^{m+a-k}+ \br{t\pm r}^{m+a-k} \right] \\
&\leq C\e r^{-m-a} \br{t-r}^{m+a-k-|\b|}.
\end{align*}
Invoking the equivalence of $r$ with $\br{t+r}$, we may thus deduce the desired estimate.

Finally, the uniqueness assertion of our theorem can be established using the same energy argument as in \cite{KKo}.
Namely, the condition that \cite{KKo} imposed was of the form
\begin{equation}\label{kko}
|\d_r u_0 (r,t)| + |\d_t u_0(r,t)| = O(r^{-(n-1)/2+\de}) \quad\quad\text{as\, $r\to 0$}
\end{equation}
for some $\de>0$.  In our case, the estimates we just proved are such that
\begin{equation*}
|\d_r u_0 (r,t)| + |\d_t u_0(r,t)| = O(r^{l-1-m}) \quad\quad\text{as\, $r\to 0$.}
\end{equation*}
However, $l\geq 1$ and $m$ is either $(n-3)/2$ or $(n-2)/2$, so the condition \eqref{kko} does hold.
\end{proof_of}

\begin{corollary}\label{imp}
Using the assumptions and notation of Theorem \ref{dec}, one may improve the estimates provided for the decay rates
$0\leq k< m+a= (n-1)/2$ of part (a) as follows.
\begin{itemize}
\item[(i)]
When $0\leq k< m+a-l$, we actually have
\begin{equation*}
|D^\b u_0(r,t)| \leq C_0\e r^{l-|\b|-m} \cdot \br{t+r}^{m-l-k}.
\end{equation*}

\item[(ii)]
When $k_0-1 \leq k-m-a+l< k_0$ for some integer $1\leq k_0\leq l$, we have
\begin{equation*}
|D^\b u_0(r,t)| \leq C_0\e r^{l-|\b|-m} \cdot \br{t-r}^{\min(l-|\b|-k_0,0)} \br{t+r}^{m-l-k-\min(l-|\b|-k_0,0)}.
\end{equation*}

\end{itemize}
\end{corollary}

\begin{proof}
If either $t\geq 2r$ or $r\leq 1$, then $\br{t+r}$ is equivalent to $\br{t-r}$ and each of the desired estimates follows
from part (a) of Theorem \ref{dec}.  Let us now assume that $t\leq 2r$ and $r\geq 1$.  As in Case 2 of the previous
proof, our task reduces to the estimation of
\begin{equation}\label{B13}
\mathbb{B}_1 + \mathbb{B}_2 = C_1'\e r^{-m-a} \left[ \int_{|t-r|}^{t+r} \frac{\la^{l-|\b|} \br{\la}^{m-l-k} d\la}{(r-t
+\la)^{1-a}} + |t\pm r|^{a+l-|\b|} \br{t\pm r}^{m-l-k} \right].
\end{equation}
In what follows, we only concern ourselves with part (ii) because part (i) is easier to settle.  Given any $|t-r|\leq
\la\leq t+r$ and any $0\leq k_0\leq l$, let us write
\begin{equation*}
\la^{l-|\b|} = \la^{\min(l-|\b|, k_0)} \cdot \la^{\max(l-|\b|-k_0,0)}.
\end{equation*}
Since $|\b|\leq l$ by assumption, the exponents on the right hand side are non-negative, so
\begin{align*}
\la^{l-|\b|} &\leq \br{\la}^{\min(l-|\b|-k_0,0)+k_0} \cdot (t+r)^{\max(l-|\b|-k_0,0)} \\
&\leq \br{\la}^{k_0} \cdot \br{t-r}^{\min(l-|\b|-k_0,0)} \cdot (3r)^{\max(l-|\b|-k_0,0)}
\end{align*}
whenever $|t-r|\leq \la\leq t+r$ and $t\leq 2r$.  Inserting this fact in \eqref{B13}, we thus arrive at
\begin{align*}
\mathbb{B}_1 + \mathbb{B}_2 &\leq C\e r^{-m-a+\max(l-|\b|-k_0,0)} \cdot \br{t-r}^{\min(l-|\b|-k_0,0)} \times\\
&\quad\quad \left[ \int_{|t-r|}^{t+r} \frac{\br{\la}^{m-l-k+k_0} d\la}{(r-t +\la)^{1-a}} + \br{t\pm r}^{m-l-k+k_0+a}
\right].
\end{align*}
For the decay rates of part (ii), we have $m-l-k+k_0>-a$, so Lemma \ref{ts1} ensures that
\begin{equation*}
\int_{|t-r|}^{t+r} \frac{\br{\la}^{m-l-k+k_0} d\la}{(r-t +\la)^{1-a}} + \br{t\pm r}^{m-l-k+k_0+a} \leq C\br{t+r}^{m-l-k
+k_0 +a}.
\end{equation*}
Since $r$ is equivalent to $\br{t+r}$ when $r\geq \max(t/2,1)$, the desired estimate follows easily.
\end{proof}

\section*{Acknowledgement}
The author would like to express his gratitude to Walter A.~Strauss, under whose direction this research was conducted as
part of the author's doctoral dissertation at Brown University.

\end{document}